\documentclass[12pt]{article}
\usepackage{amsmath}
\usepackage{amssymb}
\usepackage{amscd}
\usepackage{amsthm}
\usepackage{mysects}

\def\no{\noindent}

\numberwithin{equation}{subsection}

\theoremstyle{plain}
\newtheorem{theorem}[equation]{Theorem}

\newtheorem{proposition}[equation]{Proposition}
\newtheorem{lemma}[equation]{Lemma}
\newtheorem{corollary}[equation]{Corollary}

\theoremstyle{remark}

\theoremstyle{definition}

\newtheorem{question}[equation]{Question}

\newtheorem{conjecture}[equation]{Conjecture}

\newcommand{\ra}{\rightarrow}
\newcommand{\restr}{\mbox{\Large \(|\)\normalsize}}
\newcommand{\N}{\mathbb N}

\newcommand{\R}{\mathbb R}

\newcommand{\Z}{\mathbb Z}

\renewcommand{\H}{\mathbb H}

\newcommand{\acts}{\curvearrowright}

\newcommand{\al}{\alpha}

\def\De{\Delta}

\def\eps{\epsilon}
\def\ga{\gamma}
\def\Ga{\Gamma}

\def\la{\lambda}

\def\lang{\langle}

\def\for{\quad\mathrm{for}\quad}
\def\diam{\mathrm{diam}}
\def\dist{\mathrm{dist}}

\def\length{\mathrm{length}}
\def\:{\colon}
\def\sub{\subset}
\def\sph{{\mathbb S}^2}
\def\St{\mathrm{St}}
\def\mesh{\mathrm{mesh}}
\def\Mod{\mathrm{Mod}}
\def\GMod{\hbox{\rm mod}}

\def\Om{\Omega}

\def\ra{\rightarrow}
\def\Ra{\Rightarrow}
\def\rang{\rangle}

\def\geo{\partial_{\infty}}

\def\defeq{:=}

\newcommand{\M}{M\"obius}
\newcommand{\qm}{quasi-M\"obius}
\newcommand{\qs}{quasisymmetric}

\begin{document}
\nocite{*}

\title{Quasisymmetric parametrizations of two-dimensional metric
 spheres}
\author{Mario Bonk\thanks{Supported by a Heisenberg fellowship
of the Deutsche Forschungsgemeinschaft.}\ \ 
 and Bruce Kleiner\thanks{Supported by NSF grant DMS-9972047.} }
\date{July 16, 2001}
\maketitle

\subsection{Introduction}

According to the classical uniformization theorem, 
every smooth Riemannian surface $Z$ homeomorphic to the  $2$-sphere  
is conformally diffeomorphic to  $\sph$ 
(the unit sphere in $\R^3$ equipped with the Riemannian metric 
induced by the ambient Euclidean metric). 
The availability of a  similar uniformization procedure  is 
highly desirable in a nonsmooth setting, in particular in connection with
Thurston's hyperbolization conjecture.
This question of nonsmooth uniformization 
was addressed by Cannon in his combinatorial Riemann mapping 
theorem \cite{cannonacta}.
He considers topological surfaces equipped with 
combinatorial data that lead to a notion of approximate 
conformal moduli of rings. He then finds conditions on the combinatorial
structure that imply the existence of coordinates systems
 on the surface that 
relate these combinatorial moduli to classical analytic  moduli
in the plane.

In this paper we present a different approach to  nonsmooth
uniformization. 
We start with a metric space $Z$ homeomorphic to $\sph$ and ask  
for conditions under which $Z$ can be mapped onto $\sph$ by a quasisymmetric
homeomorphism. The class of quasisymmetries is an appropriate analog 
of conformal\footnote{A homeomorphism between compact Riemannian manifolds
is quasisymmetric iff it is quasiconformal.  There does not seem to be any hope of
a general
 existence theory for {\em conformal} mappings beyond the Riemannian
setting: by any reasonable definition, two norms on $\R^2$ define
locally conformally equivalent metrics iff the  norms are 
isometric.}
mappings in a metric space context.
Quasisymmetric homeomorphisms also arise in the theory
of Gromov hyperbolic metric spaces -- quasi-isometries
between Gromov hyperbolic spaces induce quasisymmetric boundary homeomorphisms.  
Our  setup has the advantage that
we can exploit recent  notions and methods from the analysis
on metric spaces.  
Our main result Theorem \ref{suff}  gives a necessary and sufficient for $Z$ 
to be quasisymmetric equivalent to $\sph$. Since the formulation of this 
theorem requires some preparation, we  postpone this to
Section \ref{asymptoticconditions}.
 In this introduction we  formulate two consequences  of our methods
 that are  easier to state. The first 
result
 answers a question of Heinonen  and Semmes affirmatively
(cf.\ \cite{heinsem}, Question 3 and  \cite{sem4}, Section 8)
and was the original motivation for this
paper.
 
\begin{theorem}
\label{mainthm}
Let $Z$ be an Ahlfors $2$-regular metric space homeomorphic
to $\sph$.
Then $Z$ is quasisymmetric to $\sph$ if and only if
$Z$ is linearly locally contractible.
\end{theorem}
 
\no
We recall that a metric space $X$ is Ahlfors $Q$-regular
if there is a constant $C>0$ such that the $Q$-dimensional
Hausdorff measure  ${\cal H}^Q$ of every open $r$-ball $B(a,r)$
satisfies
 
$$C^{-1}r^Q\leq {\cal H}^Q(B(a,r))\leq Cr^Q,$$
when $0<r\le\diam(X)$.
 A metric space is linearly locally
contractible if there is a constant $C$ such that
every  small ball is contractible inside
a ball whose radius is $C$ times larger; for closed
surfaces linear local contractibility is equivalent
to linear local connectedness, see Section \ref{cross-ratios}.

The statement of  Theorem \ref{mainthm} is quantitative  in a sense 
that will be explained below (See the comment after the proof of Theorem 
\ref{mainthm} in Section \ref{theproofsofthetheorems}).

The  problem considered here  is just a special case of the general problem of
characterizing a metric space $Z$ up to quasisymmetry.
Particularly interesting are the cases when $Z$ is  
 $\R^n$ or the standard sphere ${\mathbb S}^n$.
Quasisymmetric characterizations of $\R$ and ${\mathbb S}^1$ have been given
by Tukia and V\"ais\"al\"a \cite{tukvai}. If $n\ge 3$ then results by Semmes 
 \cite{sem2} show that natural conditions which one might expect to
imply that a metric space is \qs\ to   ${\mathbb S}^n$ 
(or $\R^n$), are in fact insufficient;  at present these cases
look intractable.

A result similar to Theorem \ref{mainthm}
has been proved by Semmes \cite{sem3}
 under the additional assumption  that $Z$ is  a smooth Riemannian surface.  
The hypothesis of $2$-regularity in the  theorem
is essential.  A metric $2$-sphere  containing an open set   bilipschitz
equivalent to the unit disk $B(0,1)\subset \R^2$ with the metric

$$d_\al((x_1,y_1),(x_2,y_2))=|x_1-x_2|+|y_1-y_2|^\al,$$
where  $0<\al<1$, will never be \qs ally homeomorphic to $\sph$,
see \cite{tuk2,vai2}.
We also mention  that the construction of Laakso \cite{laa} provides examples
of Ahlfors $2$-regular, linearly locally contractible  $2$-spheres
which are not bilipschitz homeomorphic to $\sph$; this 
shows that one cannot replace the word ``\qs '' with 
``bilipschitz'' in the statement of the theorem.
Finally we point out that the $n$-dimensional analog
of Theorem \ref{mainthm} is false for $n>2$ according to the results
by Semmes  \cite{sem2}:   for $n>2$ there are 
 linearly locally contractible and $n$-regular metric $n$-spheres
 which are not \qs\ to the standard $n$-sphere.
However, if an $n$-regular  $n$-sphere
admits an appropriately large group of symmetries, then
it must be \qs ally homeomorphic to the standard $n$-sphere,
see \cite{small}. 

Theorem \ref{mainthm}  is closely related to a theorem of Semmes \cite{sem}
which  shows that 
an  Ahlfors $n$-regular metric space  that is a linearly locally 
contractible  topological $n$-manifold  
satisfies a $(1,1)$-Poincar\'e inequality (see Section \ref{modulus})  and 
hence  has nice analytic properties. 
His result shows 
in particular that a $2$-sphere as in our theorem satisfies a
Poincar\'e inequality. 
We will not use this result, since it does not substantially simplify
our arguments, and in
fact our theorem together with a result by Tyson \cite{tys}
gives a different  way to establish  a Poincar\'e inequality in our case.
Our methods could also easily be adapted to show this directly. 

From an analytic  perspective it is interesting
to consider metric spaces that  
satisfy Poincar\'e inequalities by assumption (cf.\ 
\cite{heinkosk,sem,CHEEGER,bp1,bp2,laaksopi}). For an 
 Ahlfors $Q$-regular metric
space a  $(1,Q)$-Poincar\'e inequality is equivalent to
the $Q$-Loewner property as introduced by Heinonen and
Koskela \cite{heinkosk}, see  Section \ref{modulus}.  It turns
out that in dimension $2$, this is a very restrictive condition:

\begin{theorem}
\label{thmloew}
Let $Q\ge 2$ and  $Z$ be  an Ahlfors $Q$-regular metric space homeomorphic
to $\sph$. If $Z$ is  $Q$-Loewner, 
then
$Q=2$ and $Z$ is quasisymmetric to $\sph$.
\end{theorem}

\no
By the result of Semmes \cite{sem} the space $Z$ will actually
satisfy a $(1,1)$-Poincar\'e inequality. 

The analog of Theorem \ref{thmloew} in higher dimensions
is false---one has the examples
of Semmes cited above.  Also, the standard Carnot metric on
the $3$-sphere is Ahlfors $4$-regular and $4$-Loewner. 
In  view of these examples one can summarize Theorem \ref{thmloew} by 
saying that there are no  exotic geometric structures 
on $\sph$ that are analytically nice.  

Another source of examples of Ahlfors regular, linearly locally contractible
metric spheres is the theory of Gromov hyperbolic groups.  
The boundary $\geo G$
of a hyperbolic group $G$ has a natural family of Ahlfors regular metrics
which are \qs\ to one another by the identity homeomorphism.
When $\geo G$ is homeomorphic to a sphere, then these metrics
are all linearly locally contractible.  Cannon \cite{cannonacta} has
conjectured that when $\geo G$ is homeomorphic to $\sph$,
then $G$ admits a discrete, cocompact, and isometric action
on hyperbolic $3$-space $\H^3$.   This conjecture is a major
piece of Thurston's hyperbolization conjecture for 
$3$-manifolds\footnote{The Hyperbolization Conjecture is
part of the full Geometrization Conjecture.  It says
that a closed, irreducible, aspherical $3$-manifold 
admits a hyperbolic structure provided its fundamental
group does not contain a copy of $\Z\times \Z$.}.
By a theorem of Sullivan \cite{sul} 
Cannon's conjecture is equivalent to the following conjecture:
\begin{conjecture}
If $G$ is a hyperbolic group and $\geo G$ is homeomorphic to $\sph$, 
then $\geo G$ (equipped with one of the metrics mentioned above) 
is \qs\ to $\sph$.
\end{conjecture}
In this connection we raise the following question:

\begin{question}
\label{loewnerq}
Suppose $G$ is an infinite hyperbolic group, and neither $G$ nor any
finite index subgroup of $G$ splits over a virtually 
cyclic group. 
 Is $\geo G$ quasi-symmetric
to a $Q$-regular metric space which satisfies a
$(1,Q)$-Poincar\'e inequality, for some $Q$?
\end{question}

\no
Note that by  work of Bestvina-Mess, Bowditch, and Swarup,
a Gromov hyperbolic group $G$  does not virtually split over a virtually
cyclic group if and only if 
$G$ is  non-elementary
and has a  connected boundary with no local cut points.
An affirmative answer to  Question \ref{loewnerq} would imply Cannon's
conjecture, by Sullivan's theorem and Theorem \ref{thmloew}.

We now turn to the problem of finding necessary and sufficient
conditions for a metric space to be \qs\ to $\sph$.
It follows easily from the definitions that a compact
metric space $X$ which is \qs\ to a doubling (respectively
linearly locally contractible) metric space is itself
doubling (respectively linearly locally contractible).
Therefore any metric space \qs\ to a standard sphere
is doubling and linearly locally contractible.  
In Section \ref{theproofsofthetheorems} we give two
different necessary and sufficient conditions for a metric
$2$-sphere to be \qs\ to $\sph$, Theorems \ref{necsuff1}
and \ref{necsuff}.   Roughly speaking, Theorem \ref{necsuff} 
says that a doubling,
linearly locally contractible metric $2$-sphere $Z$ is
\qs\ to $\sph$ if and only if, when one consider a 
sequence of finer and finer ``graph approximations'' of $Z$, 
the corresponding combinatorial moduli of any pair of 
continua $(E,F)$ are small
provided the relative distance  $\Delta(E,F)$ as defined in (\ref{reldist})
is big.  Theorem \ref{necsuff1} is similar, except that
one assumes instead that if the moduli of the pair
$(E,F)$ are  small then the
relative distance $\De(E,F)$ is big.
We refer the reader to Section \ref{theproofsofthetheorems}
for the precise statements of these two theorems.

The problem of finding necessary and sufficient
conditions for a metric sphere to be \qs\ to $\sph$ has
some features in common with
Cannon's work  \cite{cannonacta} on the combinatorial Riemann mapping
theorem. We will discuss
this in  Section \ref{asymptoticconditions}. In this section we prove
Theorem \ref{suff} which is an improvement of Theorem \ref{necsuff}. 
One can use Theorem \ref{suff} to 
verify that certain self-similar examples are \qs\ to
$\sph$. 

We now outline   the proofs of Theorems \ref{mainthm} 
and \ref{thmloew}.  

 The first step is to use the linear local 
contractibility to produce  an embedded graph with 
controlled geometry which
 approximates  our space $Z$ on a given 
scale. This can actually be done for any
doubling, linearly locally connected metric space. If $Z$ is a 
topological $2$-sphere, then 
we can obtain a graph approximation which is, in addition, the $1$-skeleton 
of a triangulation.  In the second step we apply
a uniformization procedure.  We 
 invoke the  circle packing theorem of
Andreev-Koebe-Thurston, 
which ensures that every 
triangulation of the $2$-sphere is combinatorially equivalent to the 
triangulation
dual to some circle packing, and then map each vertex
of the graph to the center of the associated circle.
In this way we get a  mapping $f$ from the vertex
set of our approximating graph to the sphere.\footnote{Alternatively,
one  can  use  the classical uniformization
theorem to produce such a map.
To do this, one endows the sphere with a piecewise flat metric so that each
$2$-simplex of the topological triangulation is isometric to an equilateral
Euclidean triangle with side length $1$.  Such a piecewise flat
metric defines a flat Riemannian surface with isolated conical singularities,
and one can then apply the classical uniformization theorem to 
get a map  from this Riemann surface to $\sph$.}
The way to think about the 
map is that it provides
a coarse conformal change of the metric:
the scale attached to a given vertex of the graph approximation is changed 
to the scale 
given by the radius of the corresponding disc in the circle packing. 
The third step  is to show that (after suitably
normalizing the circle packing) the mapping  $f$ has
controlled \qs\ distortion.  Since in some sense $f$  
changes the metric conformally, we control its \qs\ distortion 
(in fact it is the \qm\ distortion which enters more
naturally) via modulus estimates.  There are two main
ingredients in our implementation of this idea---the 
Ferrand cross-ratio (cf.\ \cite{fer,bp2}), which mediates between the \qs\ 
distortion and the ``conformal'' distortion, and  a modulus
comparison proposition which allows one to relate (under
suitable conditions) the 
$2$-modulus of a pair of continua $E,F\subset Z$
with the combinatorial $2$-modulus of their discrete approximations 
in the approximating graph.  In the final step we take
a sequence of  graph approximations 
at finer and finer scales, and  apply 
Arzela-Ascoli to see that the corresponding mappings subconverge to a \qs\
homeomorphism from $Z$ to $\sph$.

We suggest that readers who are unfamiliar with modulus arguments
read Sections \ref{cross-ratios}, \ref{quasi-mobiusmaps}, \ref{modulus}, and Proposition \ref{FCR0}.
The proposition is a simplified version of later arguments
which bound \qm\ distortion.

\tableofcontents

\subsection{Cross-ratios}
\label{cross-ratios}
 
\label{qm}
We use the notation $\N=\{1,2,3,\dots\}$, $\N_0=\{0,1,2,\dots\}$,
 $\R^+=(0,\infty)$, and $\R^+_0=[0,\infty)$.

Let $(Z,d)$ be a metric space.  We denote by $B_Z(a,r)$ and by $\bar B_Z(a,r)$
the open and closed ball in $Z$ centered at $a\in Z$ of radius $r>0$,
respectively. 
We drop the subscript $Z$ if the space $Z$ is understood.  

The {\em cross-ratio},
$[z_1,z_2,z_3,z_4]$,
of a four-tuple of  distinct points $(z_1,z_2,z_3,z_4)$ in $Z$
is the quantity
$$[z_1,z_2,z_3,z_4]:=\frac{d(z_1,z_3)d(z_2,z_4)}{d(z_1,z_4)d(z_2,z_3)}.$$

Note that
\begin{equation} \label{cr1}
[z_1,z_2,z_3,z_4]= [z_2,z_1, z_3,z_4]^{-1}=[z_1,z_2,z_4,z_3]^{-1}=
[z_3,z_4, z_1,z_2]. 
\end{equation}
 
It is convenient to have a  quantity that is quantiatively equivalent to the
cross-ratio and has a geometrically  more
transparent meaning.
Let $a\vee b:=\max\{a,b\}$ and $a\wedge b:=\min\{a,b\}$ for
$a,b\in \R$. If $(z_1,z_2,z_3,z_4)$ is a four-tuple of distinct points in
$Z$ define
\begin{equation} \label{cr1a} \lang z_1,z_2,z_3,z_4\rang
:= \frac{d(z_1,z_3)\wedge d(z_2,z_4)}{d(z_1,z_4)\wedge
d(z_2,z_3)}.
\end{equation}
 
Then the following is true.
\begin{lemma}
\label{cr2} Let $(Z,d)$ be a metric space and
 $\eta_0(t)=3(t\vee \sqrt t)$ for $t> 0$.
Then for every four-tuple $(z_1,z_2,z_3,z_4)$ of distinct points in
$Z$ we have
\begin{equation}  \label{cr3}
\lang z_1,z_2,z_3,z_4\rang  \le \eta_0([z_1,z_2,z_3,z_4]).
\end{equation}
\end{lemma}
 
\proof  Suppose there is a four-tuple $(z_1,z_2,z_3,z_4)$
 for which the left hand side in (\ref{cr3})  is bigger
than the  right hand side. Let $t_0=[z_1,z_2,z_3,z_4]$.
 We may assume $d(z_1,z_3)\le d(z_2,z_4)$.
Then
\begin{eqnarray*}
d(z_1,z_4) & \le & d(z_1,z_3) + d(z_3, z_2) + d(z_2,z_4) \\
           & \le & 2d(z_2,z_4)+ d(z_2,z_3).
\end{eqnarray*}
Similarly, $ d(z_2,z_3)\le 2d(z_2,z_4) + d(z_1,z_4)$ and so
\begin{eqnarray*}
 d(z_1,z_4) \vee d(z_2,z_3) &\le&
 2 d(z_2,z_4) + d(z_1,z_4) \wedge d(z_2,z_3)\\
&\le & \left(2+\frac 1{\eta_0(t_0)}\right)  d(z_2,z_4).
\end{eqnarray*}
Hence,
$$ t_0=[z_1,z_2,z_3,z_4]\ge 
\frac {d(z_1,z_3)\eta_0(t_0)} {(d(z_1,z_4)\wedge d(z_2,z_3))( 1+2 \eta_0(t_0))}
\ge \frac{ \eta_0(t_0)^2}
{ 1+2 \eta_0(t_0)}> t_0.$$
This is a contradiction. \qed
 
\medskip
Using the symmetry properties (\ref{cr1}) for the cross-ratio which are also
true for the modified cross-ratio defined in (\ref{cr1a}), we obtain an
inequality as in  (\ref{cr3}) with the roles of the cross-ratios reversed
and the function $\eta_0$ replaced by the function
 $t\mapsto 1/\eta_0^{-1}(1/t)$. In particular, we conclude
that $[z_1,z_2,z_3,z_4]$ is small if and only if
$\lang z_1,z_2,z_3,z_4\rang $ is small, where the quantitative dependence
is given by universal functions. 

A metric space $(Z,d)$ 
 is called $\la$-{\em linearly locally contractible} where $\la\ge 1$,
 if 
every ball $B(a,r)$ in $Z$  with $0<r\le \diam(Z)/\la$ 
is contractible inside $B(a,\la r)$,  i.e., there exists a continuous
map $H\: B(a,r)\times [0,1] \ra B(a,\la r)$ such that $H(\cdot, 0)$ is the 
identity on  $B(a,r)$ and $H(\cdot, 1)$ is a constant map. 
The space is called linearly locally contractible, if it is 
$\la$-linearly locally contractible for some $\la\ge 1$. Similar language
will be employed for other notions that depend on numerical parameters. 

A metric space $(Z,d)$ is called $\la$-$LLC$ for $\la\ge 1$ ($LLC$ stands for
  linearly locally connected) 
if the following two conditions are satisfied:

\smallskip \noindent 
($\la$-$LLC_1$)$\quad$
 If $B(a,r)$ is a ball in $Z$ and $x,y\in B(a,r)$, then there 
exists a continuum $E\sub B(a,\la r)$ containing $x$ and $y$. 

\smallskip \noindent 
($\la$-$LLC_2$)$\quad$
 If $B(a,r)$ is a ball in $Z$ and $x,y\in Z \setminus B(a,r)$,
then there exists a con\-tinuum $E\sub Z \setminus B(a,r/\la)$ containing 
$x$ and $y$.  
\smallskip

We remind the reader that a continuum is a compact connected set consisting
of more than one point.

Linearly local contractibility implies the $LLC$
condition for compact connected topological
$n$-manifolds, and is equivalent to it when $n=2$:

\begin{lemma}
Suppose $Z$ a metric space which  is a 
compact connected topological $n$-manifold. Then: 
\begin{itemize}
\item[{\rm (i)}] 
 If $Z$ is $\la$-linearly locally contractible, then $Z$
is $\la'$-$LLC$ for each $\la' >\la$.
\item[{\rm (ii)}]  If $n=2$ and $Z$ is $LLC$, then $Z$ is linearly locally
contractible.  The linear local contractibility constant
depends on $Z$ and not just on the $LLC$ constant.
\end {itemize}
\end{lemma}

\proof
(i)
We first verify the $LLC_1$ condition.   If $a\in Z$, and 
$r>\diam(Z)/\la$, then $B(a,\la r)= Z$, so in this case  the 
$\la$-$LLC_1$ condition follows from the connectedness of $Z$.
If $r\leq \diam(Z)/\la$, then the inclusion  $i\:B(a,r)\ra B(a,\la r)$ 
is homotopic to a constant map. Hence it 
induces the  zero homomorphism on reduced $0$-dimensional homology, which
means that $\la$-$LLC_1$ holds.

Let $\la'>\la$. To see that
$\la'$-$LLC_2$ holds, we have to show that if $B(a,r')\sub Z$ is a ball with
$Z\setminus B(a,r')\ne \emptyset$, then the inclusion 
$i\: Z\setminus B(a,r') \to Z\setminus B(a,r'/\la')$ induces 
the zero homomorphism 
\begin{equation} \label{lin1}
\tilde H_0(Z\setminus B(a,r'))\stackrel{0}{\ra}
\tilde H_0(Z\setminus  B(a,r'/\la'))
\end{equation} 
for reduced singular homology with coefficients in $\Z_2$. 
Note that $Z\setminus B(a,r')\ne \emptyset$ implies $r'<\diam(Z)$. 
Moreover, we can find $0<r<r'$ close enough to $r'$ such that 
$\bar B(a,r'/\la') \sub B(a, r/\la)$.
Let $K_1:= \bar B(a,r'/\la')$ and 
$K_2:=\bar B(a,r)$. Then $K_1$ and $K_2$ are compact, and we 
have  $B(a, r'/\la')\sub K_1\sub K_2 \sub B(a,r')$.
So in  order to show (\ref{lin1}), it is enough to show
that the inclusion $i\: 
Z\setminus K_2\ra Z\setminus K_1$ induces the zero homomorphism
\begin{equation} \label{lin2}
\tilde H_0(Z\setminus K_2)\stackrel{0}{\ra}
\tilde H_0(Z\setminus K_1). 
\end{equation} 
It follows from the path connectedness of $Z$ and the long exact sequence  
for singular homology that the natural map $\partial\: H_1(Z, Z\setminus K_i)
\ra \tilde H_0(Z\setminus K_i) $ is surjective for $i\in \{1,2\}$. 
Hence (\ref{lin2}) is true, if the inclusion $i\: 
(Z, Z\setminus K_2)\ra (Z, Z\setminus K_1)$ induces the zero homomorphism 
\begin{equation} \label{lin3}
 H_1(Z, Z\setminus K_2)\stackrel{0}{\ra}
 H_1(Z, Z\setminus K_1). 
\end{equation} 
Now duality \cite[Theorem 17, p.\ 296]{spanier} shows that  for each compact
subset $K\sub Z$ we have an isomorphism
$H_1(Z,Z-K)\simeq\check{H}^{n-1}(K)$, where $\check{H}^{*}$ denotes \v{C}ech
cohomology with coefficients in $\Z_2$. This isomorphism is natural, and 
hence compatible with inclusions.   Since 
$K_1 \sub  B(a, r/\la) \sub  B(a, r) \sub K_2$ and $r<r'< \diam(Z)$,
it follows from our assumptions that $K_1$ contracts to 
a point inside $K_2$. Hence the inclusion $i:K_1\ra K_2$ 
induces the zero homomorphism $\check{H}^{n-1}( K_2)\stackrel{0}{\ra}
 \check{H}^{n-1}( K_1).$ Therefore, (\ref{lin3}) holds which implies
(\ref{lin1}) as we have seen.  

(ii) It is enough to show that the inclusion $i\: B(a,r) \ra B(a, \la r)$ 
is homotopic to a constant map, if $r>0$ is sufficiently small 
independent of $a\in Z$.
 Since $Z$ is a compact $2$-manifold, every sufficiently small ball
lies precompactly in an open subset of $Z$ homeomorphic to  $\R^2$.
So without loss of generality we may assume that the sets $U:= B(a,r)$
and $V:=B(a, \la r)$ are bounded and open  subsets of $\R^2$ with $U\sub V$. 
Now  $\la$-$LLC_1$  implies that $U$ lies in a single component 
of $V$. So in order to show that  $U$ is contractible inside 
$V$, it is enough to show that each component $\Omega$ of $U$ is contained in 
a simply connected (and hence contractible) subregion of $V$.

 The condition $\la$-$LLC_2$ implies that $\R^2\setminus V$ lies in one,
namely the unbounded component of $\R^2\setminus U$. 
It follows in particular that if $\ga$ is a Jordan curve in $U$, then 
the interior region  $I(\ga)$ of $\ga$  is contained in $V$. 

A well-known fact from plane topology is that every region $\Omega$ can be
written as an nondecreasing  union $\Om= \bigcup_{i=0}^\infty  \Omega_i$, 
where $\Omega_i$ is a region with $\bar \Omega_i \sub\Om$
whose boundary consists 
of finitely many Jordan curves.
One of the boundary components $\ga_i$ of $\Om_i$ is a Jordan curve
whose interior $I(\ga_i)$ contains $\Om_i$.
Now if $\Om$ is a component of $U$, then $\ga_i \sub \Om \sub U$, and 
so  $I(\ga_i) \subset V$ as we have seen. Hence $\Om\sub \bigcup_{i=0}^\infty 
I(\ga_i)\sub V$
lies in the union  of a  nondecreasing sequence of Jordan subregions
of $V$. This union is a simply connected subregion of $V$ containing 
$\Om$. \qed

In view of the lemma we prefer to
work with the weaker $LLC$-condition instead of
linear local contractibility in the following.

If $E$ and $F$ are  
continua in $Z$  we denote by
\begin{equation} \label{reldist}
   \Delta(E,F) \defeq {\dist(E,F) \over \diam(E) \wedge
\diam(F)}
\end{equation}
the {\em relative distance} of $E$ and $F$.
 
\begin{lemma}\label{goodcont} Suppose $(Z,d)$ is $\la$-$LLC$. 
Then there exist  functions $\delta_1,\delta_2\:\R^+ \ra \R^+$  
depending only on $\la$ 
with the following properties. Suppose $\eps>0$ and $(z_1,z_2,z_3,z_4)$ is
a four-tuple of distinct points in $Z$.  

\begin{itemize}

\item[{\rm (i)}] If $[z_1,z_2,z_3,z_4] <\delta_1(\eps)$, then there exist 
continua $E,F\sub Z$ with $z_1, z_3\in E$, $z_2, z_4\in F$ and 
$\Delta(E,F)\ge 1/\eps$.

\item[{\rm (ii)}] If there exist continua $E,F\sub Z$ with $z_1, z_3\in E$,
$z_2, z_4\in F$ and $\Delta(E,F)\ge 1/\delta_2(\eps)$, then 
$[z_1,z_2,z_3,z_4] <\eps$.  

\end{itemize}
\end{lemma}

As the proof will show,
the function $\delta_2$ can actually be chosen as a numerical 
function independent of $\la$. 

\proof We have to show that  
$[z_1,z_2,z_3,z_4]$ is small, if and only if  there exist two 
continua with large 
relative distance  containing 
$\{z_1, z_3\}$ and $\{z_2, z_4\}$, respectively.  
 
Suppose $s=[z_1,z_2,z_3,z_4]$ is small. Then by Lemma  \ref{cr2}
the quantity
\begin{equation}\label{FRCL20}
t:=\lang z_1,z_2,z_3,z_4\rang
= \frac{d(z_1,z_3)\wedge d(z_2,z_4)}{d(z_1,z_4)\wedge
d(z_2,z_3)}.
 \end {equation}
is small,  quantitatively. We may assume $t<1$ and 
 $r:= d(z_1, z_3)\le d(z_2, z_4)$.
Since $Z$ is $\la$-{\em LLC} and $z_1,z_3\in B(z_1, 2r)$, there
exists a continuum $E$ connecting $z_1$ and $z_3$ in $B(z_1, 2\la r)$.
Let $R:=r(1/t-1) >0$. Then $d(z_1, z_4) \ge r/t> R$ and 
$ d(z_1, z_2) \ge d(z_2, z_3)-d(z_1,z_3)\ge r(1/t-1)=R.$ 
Thus  $z_2,z_4$ are
in the complement of $B(z_1,R)$, and so  there exists a continuum $F$
 connecting $z_2$ and $z_4$ in $Z\setminus B(z_1, R/\la )$.
For the relative distance of $E$ and $F$ we get
$$ \Delta(E,F) = \frac{\dist(E,F)}{\diam(E) \wedge \diam(F)} \ge
\frac {R/\la-2\la r}{4\la r} > 1/(4\la^2t)-1,
$$
which is uniformly large if $s$ and so $t$  are small. 

Now suppose that there exist continua $E,F\sub Z$ with 
 with $z_1,z_3\in E$ and $z_2,z_4\in F$ for which
$\Delta(E,F)$ is large.
Since 
$$ \lang z_1,z_2,z_3,z_4\rang
= \frac{d(z_1,z_3)\wedge d(z_2,z_4)}{d(z_1,z_4)\wedge
d(z_2,z_3)} \le   \frac {\diam(E) \wedge \diam(F)}{\dist(E,F)}=
1/\Delta(E,F),
$$
we conclude from Lemma \ref{cr2} that $[z_1,z_2,z_3,z_4]$ is uniformly
small. \qed 

\medskip
In the proof of this lemma we used for the first time the expression 
``If $A$ is small, then $B$ is small, quantitatively." This and similar 
language  will be very 
convenient in the following, but it requires some explantion.
By this expression
we mean that  an inequality $B \le \Psi(A)$ for the quantities $A$
and $B$  holds, where $\Psi$ is a positive function  with $\Psi(t)\to 0$ 
if $t\to 0$ that depends   only on the {\em data}.   
The   data are  some ambient  parameters associated to the
given space, function, etc.\   In the proof above the data consisted 
just of the parameter $\la$ in the $LLC$-condition for $Z$.

\subsection{Quasi-M\"obius maps}
\label{quasi-mobiusmaps} 

Let $\eta\:\R_0^+\ra\R_0^+$
be a homeomorphism, i.e., a strictly increasing nonnegative function with 
$\eta(0)=0$,  and let   $f\:X\ra Y$
be an injective map between metric spaces.
The map $f$ is an {\em $\eta$-quasi-\M\ map} if for
every four-tuple $(x_1,x_2,x_3,x_4)$ of distinct
points in $X$, we have
$$[f(x_1),f(x_2),f(x_3),f(x_4)]\leq \eta([x_1,x_2,x_3,x_4]).$$
Note that by exchanging the roles of $x_1$ and $x_2$, one
gets the lower bound
$$\eta([x_1,x_2,x_3,x_4]^{-1})^{-1}\leq [f(x_1),f(x_2),f(x_3),f(x_4)].$$
Hence
 the inverse $f^{-1}\:f(X)\ra X$ is also \qm.

Another way to express that $f$ is \qm\ is to say that the cross-ratio
$[x_1,x_2,x_3,x_4]$ of a four-tuple of distinct points  is small
if and only if the cross-ratio
$[f(x_1), f(x_2), f(x_3), f(x_4)]$ is small, quantitatively. 
This is easy to see using the symmetry properties (\ref{cr1}) of cross-ratios.

The map $f$ is  {\em $\eta$-\qs\ } if
$$\frac{d(f(x_1),f(x_2))}{d(f(x_1),f(x_3))}
\leq \eta\left (\frac{d(x_1,x_2)}{d(x_1,x_3)}\right )$$
for every triple $(x_1,x_2,x_3)$ of distinct points in $X$.
Again it is easy to see that the inverse map $f^{-1}\: f(X)\ra X$ is 
also \qs. Two metric spaces $X$ and $Y$ are called {\em quasisymmetric},
if there exists a homeomorphism $f\: X\ra Y$ that  is quasisymmetric.

Intuitively, a quasisymmetry is a map between metric spaces that 
maps  balls to roundish objects that can be trapped 
between two balls whose radius ratio is bounded by a fixed constant.
Based on this it is easy to see  the quasisymmetric invariance 
of  properties like linear local contractibility or linear 
local connectivity.

We list some properties of \qm\ and \qs\ maps (cf.\ \cite{vais}): 
 
\begin{itemize} 

\item[(1)] Quasi-M\"obius and \qs\ maps are homeomorphisms onto their image.   

\item[(2)] The composition of an $\eta_1$-\qm\ map with an
$\eta_2$-\qm\ map is an $\eta_2\circ\eta_1$-\qm\ map.
 
\item[(3)] An $\eta$-quasisymmetric map is $\tilde \eta$-\qm\ with $\tilde 
\eta$ depending only on $\eta$.
 
Conversely, every \qm\ map between bounded spaces is quasisymmetric. 
This statement is not quantitative in general, but we have:  

\item[(4)] Suppose $(X,d_X)$ and $(Y,d_Y)$ are bounded metric spaces, 
     $f\: X\ra Y$ is $\eta$-\qm, and $\la \ge 1$.   Suppose $(x_1,x_2,x_3)$
and $(y_1,y_2,y_3)$ are triples of distinct points in $X$ and $Y$, respectively,
such that $f(x_i)=y_i$ for $i\in \{1,2,3\}$, $d_X(x_i, x_j)\ge \diam(X)/\la$
and $d_Y(y_i, y_j)\ge \diam(Y)/\la $ for $i,j\in \{1,2,3\}$, $i\ne j$. 
Then $f$ is $\tilde  \eta$-quasisymmetric with $\tilde  \eta$ 
depending only on $\eta$ and $\la$.

\end{itemize}

\medskip
We will need the following convergence 
property of \qm\ maps which we  state
as a separate lemma.

\begin{lemma}
\label{TBA}
Suppose $(X, d_X)$ and $(Y,d_Y)$ are  compact 
 metric spaces, and let
$f_k\:D_k\ra Y$ for $k\in \N$ be an $\eta$-\qm\ map defined on 
a subset $D_k$ of $X$. Suppose 
$$ \lim_{k\to \infty}\,  \sup_{x\in X}\,  \dist(x, D_k) = 0$$ 
and that  for $k\in \N$ there exists triples $(x_1^k, x_2^k,x_3^k)$
and $(y_1^k, y_2^k,y_3^k)$ of points in $D_k\sub X$
and $Y$, respectively, such that 
$f(x_i^k)=y_i^k$, $k\in \N$, $i\in \{1,2,3\}$, 
$$ \inf \{ d_X(x_i^k,x^k_j): k\in \N, \, i,j\in \{1,2,3\},\, i\ne j\} >0$$
and 
$$ \inf \{ d_Y(y_i^k,y^k_j): k\in \N, \, i,j\in \{1,2,3\},\, i\ne j\} >0.$$
Then $f_k$ subconverges uniformly 
 to an $\eta$-\qm\ map $f\:X\ra Y$, i.e.\ there exists a 
sequence $(k_n)$ in $\N$  such that 
$$ \lim_{n\to \infty} \sup_{x\in D_{k_n}} d_Y(f(x), f_{k_n}(x)) =0.  $$ 
\end{lemma}

The assumptions imply that the functions  $f_k$ are equicontinuous (cf.\ 
\cite[Thm.\ 2.1]{vais}).  The proof of the lemma
 then follows from standard arguments, 
and we leave the details  to the reader.

\begin{lemma} 
\label{sepacont} Suppose $(X,d_X)$ and $(Y, d_Y)$ are metric spaces, and
 $f\: X\ra Y$ is an  $\eta$-\qm\ map. Then there exists a function
$\Phi\:\R^+\ra \R^+$ with $\lim_{t\to \infty}\Phi(t)=\infty$ depending
only on $\eta$ such that  the following statement holds.
 
If $E,F\sub X$ are disjoint continua, then
$$ \Delta(f(E), f(F))\ge \Psi(\Delta(E,F)). $$
\end{lemma}
If $f$ is surjective, and we apply the lemma 
to the  inverse map $f^{-1}$, we get a similar inequality
with the roles of sets and images sets reversed.
These inequalities say
that  the relative distance of two continua is large if and only if
the relative disctance of the image sets under a \qm\ map is large,
quantitatively.

Since every quasisymmetric map is also \qm, this last statement is also
true for quasisymmetric maps.
\proof Let $E':=f(E)$ and $F':=f(F)$. Then $E'$ and $F'$ are continua.
Hence there exist points $y_1\in E'$ and $y_3\in F'$ with $d_Y(y_1,y_3)
=\dist(E',F')$. Moreover, we can find points $y_4\in E'$ and
$y_2\in F'$ with $d_Y(y_1,y_4)\ge \diam(E')/2$ and $d_Y(y_2,y_3)\ge
\diam(F')/2$. Then
$$ \Delta(E',F') \ge 2 \lang y_1, y_2, y_3, y_4\rang. $$
On the other hand, if $x_i:=f^{-1}(y_i)$, then
$$ \Delta(E,F) \le \lang x_1,x_2,x_3,x_4\rang$$
by the very definition of these quantities.

Now if $\Delta(E,F)$ is large, then $\lang x_1,x_2,x_3,x_4\rang$ is at 
least as  large.
Since $f$ is $\eta$-\qm\ it follows from Lemma \ref{cr2} that
$\lang y_1, y_2, y_3, y_4 \rang$
and hence $\Delta(E',F')$ are large, quantitatively.
\qed

\medskip
A metric space $(Z,d)$ is called {\em weakly $\la$-uniformly perfect},
$\la> 1$, if for every $a\in Z$ and $0<r\le \diam(Z)$ the following 
is true:  if the ball $\bar B(a,r/\la)$ contains a point distinct from 
$a$, then $B(a,r)\setminus \bar B(a,r/\la)\ne \emptyset$. 

This condition essentially says that at each point $a\in Z$ the space
 is uniformly 
perfect in the usual sense up to the  scale for which there exist points
different from $a$.

A metric space  $(Z,d)$ is called $C$-{\em doubling}, $C\ge 1$, if
every ball of radius $r>0$ can be covered by 
$C$ balls of radius $r/2$. 

\begin{lemma} \label{weakqm} 
Suppose $(X,d_X)$ and $(Y,d_Y)$
are metric spaces, and $f\: X\ra Y$ is a bijection.   
Suppose that  $X$ is weakly $\la$-uniformly perfect, $Y$ is $C_0$-doubling,
and there exists a function $\delta_0\:\R^+\ra \R^+$ such 
that
\begin{equation} \label{wq1} 
[f(x_1), f(x_2), f(x_3), f(x_4)] < \delta_0(\eps) \Ra   
[x_1, x_2, x_3, x_4] < \eps, 
\end{equation} 
whenever $\eps>0$ and $(x_1,x_2,x_3, x_4)$ is a four-tuple of distinct points
in $X$. Then $f$ is $\eta$-\qm\ with $\eta$  depending only on $\la$, $C_0$,
and $\delta_0$. 
\end{lemma} 

As we remarked above, a bijection is \qm\ if it has the property that 
a cross-ratio of four points
 is small if and only if the cross-ratio of the image points is small,
quantitatively. The lemma says that for suitable spaces this equivalence,
which consists of  implications in two directions,
can be replaced by one of these  implications.

\proof  We have to show that for every $\eps>0$ there exists $\delta=
\delta(\eps, \la, C_0, \delta_0)>0$
such that 
\begin{equation} \label{wq2}
[x_1, x_2, x_3, x_4] < \delta \Ra [f(x_1), f(x_2), f(x_3), f(x_4)] <\eps,
\end{equation}
whenever $(x_1,x_2,x_3, x_4)$ is a four-tuple of distinct points
in $X$. 
By Lemma \ref{cr2},  for this purpose it is enough to show  the following: 
if  $\eps\in (0,1]$ and  $(x_1,x_2,x_3, x_4)$
is a four-tuple of distinct points
in $X$ with $\lang x_1, x_2, x_3, x_4\rang < \delta$ and
$\lang y_1,y_2,y_3,y_4 \rang \ge 
\eps$, where $y_i=f(x_i)$, $i\in \{1,2,3,4\}$, then we obtain 
a contradiction if $\delta$ is smaller than a positive number depending
on $\eps$, $\la$, $C_0$, and $\delta_0$.  

We may assume that $s:=d_X(x_1,x_3)\le d_X(x_2,x_4)$, and 
\begin{equation} \label{wq3} 
t:= d_Y(y_1,y_4)= \min \{ d_Y(y_i, y_j): i\in \{1,3\}, \, j\in \{2,4\}\} . 
\end{equation} 
Then 
\begin{equation} \label{wq3a}
d_Y(y_i, y_j)\ge \eps t \for i,j\in \{1,2,3,4\}, \, i\ne j.
\end{equation} 

We have that 
\begin{eqnarray*}
\diam(X) &\ge & \min \{ d_X(x_i, x_j): i\in \{1,3\}, \, i\in \{2,4\}\} \\
&\ge &  d_X(x_1,x_4)\wedge 
d_X(x_2,x_3)-d_X(x_1,x_3) \, \ge \, (1/\delta-1) s. 
\end{eqnarray*}
Since we may assume that $(1/\delta-1)> \la^2$, 
 we can choose $N\in \N$ such that 
\begin{equation} \label{wq4a}
 \la^{2N}<  (1/\delta-1) \le \la^{2N+2}. 
\end{equation}

Since $X$ is weakly $\la$-uniformly perfect, $x_3\in \bar B(x_1,s)$ and 
$\la^{2N} s< \diam(X)$,  there exist points
$z_i\in X$ for $i\in \{1,\dots, N\}$  such that 
$$ z_i\in B(x_1, \la^{2i}s ) \setminus \bar B(x_1, \la^{2i-1}s).$$ 
Then  $$\dist(z_i, \{x_1,x_3\}) \le  (\la^{2i}+1)s \for 
i\in \{1,\dots, N\}$$ and 
$$ d_X(z_i,z_j)\ge \la^{2j-2}(\la-1)s \for i,j\in \{1,\dots, N\}, \, i< j. $$ 
It follows that 
\begin{equation*}\lang z_i,u,z_j,v\rang \ge c(\la)>0 \end{equation*}
whenever $i,j\in \{1,\dots, N\}$, $ i\ne j$,
$ u\in \{x_1,x_3\}$ and  $v\in \{x_2,x_4\}$.
By our hypotheses and Lemma \ref{cr2} there exists $c_1 \in (0,1]$
depending only
on $\delta_0$ and $\la$ such that 
\begin{equation}\label{wq5} 
\lang f(z_i),u,f(z_j),v\rang \ge c_1>0, 
\end{equation}
whenever 
$i,j\in \{1,\dots, N\}$, $i\ne j$,
$u\in \{y_1,y_3\}$, and $ v\in \{y_2,y_4\}$.

We claim that 
\begin{equation}\label{wq6}d_Y(f(z_i), f(z_j))\ge   c_1\eps t/3=:c_2t
\end{equation} 
for $i,j\in \{1,\dots, N\}$, $i\ne j$. 
For otherwise, by (\ref{wq3a}) we can pick $u\in \{y_1,y_3\}$ and 
$v\in  \{y_2,y_4\}$ such that 
$$\dist(\{ f(z_i), f(z_j)\}, \{u,v\}) \ge t\eps/3$$ and we get a  
contadiction to (\ref{wq5}). 
Moreover, at most one of the points $f(z_i)$ can lie outside 
$\bar B(y_1, c_3t)$ with $c_3=1+1/c_1$. For if this were 
true for $f(z_i)$ and $f(z_j)$, $i\ne j$,
then again we get a contradiction to (\ref{wq5}) with $u=y_1$ and $v=y_4$.

The doubling property of  $Y$ now shows that the number of points 
in $\bar B(y_1, c_3t)$ which are $c_2t$-separated is bounded by a constant 
$C$ depending only on $C_0$, $c_2=c_2(\eps,\la,\delta_0)$
and $c_3=c_3(\eps,\la, \delta_0)$. 
Hence $N-1\le C$.
 By (\ref{wq4a}) this leads to a contradiction if $\delta$ is 
smaller than a constant depending  on $\eps$, $\la$, $C_0$,  and $\delta_0$. 
\qed

\subsection{Approximations of metric spaces}
\label{approximationofmetricspaces}

Suppose $G$ is a graph with vertex set $V$. We assume 
that there are no loops in $G$, i.e., no vertex is connected to itself by an 
edge,  and that  two arbitrary distinct vertices
are not  connected by more than one edge.
If $v_1, v_2\in V$
are connected by an edge or are identical we write $v_1\sim v_2$.
The combinatorial structure of the  graph is completely determined 
by the vertex set $V$ and this relation $\sim$. Hence we will write 
$G=(V, \sim)$. 

A {\em chain} is a sequence  $x_1, \dots, x_n$  of vertices  
with $x_1\sim x_2\sim \dots \sim
x_n$. It {\em connects} two subsets $A\sub V$ and 
$B\sub V$ if $x_1\in A$ and $x_n\in B$.

If $x,y\in V$ we let $k_G(x,y)$ be the {\em combinatorial distance}
of $x$ and $y$, i.e., $k_G(x,y)+1$ is the smallest
cardinality $\#M$ of a chain $M$ connecting $x$ and $y$. 
If $G$ is connected, then $(V,k_G)$ is a metric space, and we  define 
 $B_G(v,r):= \{u\in V: k_G(u,v)<r\}$ and $\bar B_G(v,r):=
\{u\in V: k_G(u,v)\le r\}$ for $v\in V$ and $r>0$. We drop the subscript 
$G$ if the graph under consideration is understood.    
The cardinality of the set $\{u\in V: k_G(u,v)=1\}$ is the 
 {\em valence} of $v\in V$. The {\em valence}
 of $G$ is the supremum of the valences
over all vertices in $G$. 

Now let  $(Z,d)$ be a metric space. We consider quadruples
${\cal A}=(G,p,r,{\cal U})$, where $G=(V,\sim)$ is a graph with vertex
set $V$, $p\: V\ra Z$, $r\: V\ra \R^+$ and ${\cal U}=\{U_v:v\in V\}$
is an open cover of  $Z$ indexed by the set $V$.
We let   $p_v=p(v)$ and $r_v=r(v)$ for $v\in V$.
Let
$$N_\eps(U):=\{z\in Z: \dist(z,U)<\eps\}$$
for $U\sub Z$ and $\eps>0$,  and
define  the $L$-{\em star} of $v\in V$ with respect to ${\cal A}$
 for $L>0$    as
$${\cal A}\hbox{-}{\rm St}_L(v):=\bigcup\{U_u: u\in V,\  k(u,v)< L\}.$$
We simply write ${\rm St}_L(v)$,  if no confusion can arise.
We call ${\cal A}$ a $K$-{\em approximation} of $Z$, $K\ge 1$,
if the following conditions are satisfied:
 
\begin{itemize}
 
\item[(1)] Every vertex of $G$ has valence at most  $K$.
 
\item[(2)] $B(p_v,  r_v)\sub U_v\sub
 B(p_v,  K r_v)$ for $v\in V$.
 
\item[(3)] If $u\sim v$ for $u,v\in V$, then $U_u\cap U_v \ne \emptyset$,
and $K^{-1} r_u \le r_v \le K r_v$. If $U_u\cap U_v \ne \emptyset$ for
$u,v\in V$, then $k(u,v)<  K$.

\item[(4)]
$N_{r_v/K}(U_v) \sub {\rm St}_K(v)$
for $v\in V$. 
 
\item[(5)] If $v\in V$, $z_1,\,z_2\in U_v$, then there is a path $\gamma$ in $Z$
connecting $z_1$ to $z_2$ so that $\gamma \sub {\rm St}_K(v)$.
 
\end{itemize}
The point $p_v$ should be thought of as a basepoint of $U_v$. By condition (2)
we can think of the number $r_v$ as the ``local scale" associated with
$v$. Condition (3)
says that the local scale only changes by a bounded factor if we
move to a neighbor of a vertex. Moreover, condition (3) says that the incidence
pattern of the cover ${\cal U}$ resembles the incidence pattern
of the vertices in $G$, quantitatively.
Condition (4) says that we can thicken up a set $U_v$ by a fixed amount
comparable to  the local scale by passing to the $K$-star of $v$.
Finally, condition (5) allows us to connect any two points in $U_v$ by a curve 
contained in the $K$-star of $v$.

We point out some immediate consequences of the Conditions (1)--(5):
 
\begin{itemize}
 
\item[(6)]  If $Z$ is connected, then
$G$ is connected; this follows from
(3).
 
\item[(7)] The multiplicity of ${\cal U}$ is bounded by a constant
$C=C(K)$:  if $U_{v_1}\cap \ldots\cap U_{v_n}\neq\emptyset$ 
then $\{v_1,\ldots,v_n\}\subset \bar B(v_1,K)$ by (3),
and $\#  B(v_1,K)\leq C=C(K)$ by (1). Similarly, it can be shown 
that for fixed $L>0$, the multiplicity of the cover $\{\St_L(v): v \in 
V\}$ is bounded by a number $C=C(K,L)$. 
 
\item[(8)] For the curve $\gamma$ in (5) we have $\diam(\ga)\le C r_v$
with $C=C(K)$; this follows from (2) and (3).
 
\end{itemize}

The {\em mesh size} of the 
$K$-approximation ${\cal A}$ is defined to be 
$$ \mesh({\cal A}):= \sup_{v\in V} r_v.$$  

The next lemma shows 
 that $K$-approximations behave well under quasisymmetric maps.
\begin{lemma}\label{qsofapp}
Suppose $(X,d_X)$ and $(Y,d_Y)$ are connected metric spaces, and
 $f\:X \ra Y $ is  an $\eta$-quasisymmetric homeomorphism.
Suppose $K\ge 1$ and  
${\cal A}=((V,\sim),p,r,{\cal U})$ is a $K$-approximation of $X$.
 Assume that 
\begin{equation}
\label{meshsmall}
\mesh({\cal A})  < \diam(X)/2.
\end{equation}
 For $v\in V$ define $p'_v:=f(p_v)$, $U'_v:=f(U_v)$ and 
\begin{equation}
\label{rprdef}
r'_v:= \inf\{ d_Y(f(x) ,p'_v): x\in X,\  d_X(x,p_v)\ge r_v\}. 
\end{equation} 
Let ${\cal U}'=\{U'_v: v\in V\}.$
Then ${\cal A}'=((V,\sim),p',r',{\cal U}')$ is a $K'$-approximation of $Y$
 with $K'$ 
depending only  on $K$ and $\eta$.    
\end{lemma}

\no
We emphasize that the underlying graphs of ${\cal A}$ and ${\cal A'}$ are the
same. 

Note that   by condition (\ref{meshsmall}) the set in (\ref{rprdef})
over which the infimum is taken is nonempty. The continuity of 
$f^{-1}$ implies that $r'_v$ is positive.
The number   $r'_v$ is roughly the diameter of $U'_v$. Up to 
 multiplicative constants, it is essentially the only possible choice for
$r'_v$. Our particular definition guarantees  $B_Y(p'_v, r'_v)\sub 
f(B_X(p_v, r_v))\sub f(U_v)= U'_v$. 

Up to this ambiguity in the choice of $r'_v$, the $K'$-approximation
 ${\cal A}'$ is  canonically determined by ${\cal A}$ and the map $f$.
In this sense we can say that ${\cal A}'$ is the ``image" of ${\cal A}$
under $f$.  

\proof We denote image points under $f$ by a prime, i.e., $x'=f(x)$
for $x\in X$. We also denote
 by $K_1, K_2, \dots$ positive 
constants that can be chosen only to depend on $\eta$ and $K$.

Since $X$ is connected  and the complement
of $B_X(p_v,r_v)$ is nonempty, there exist a  point  $x_v\in X$ with with 
$d_X(x_v,p_v)=r_v$. The  quasisymmetry of $f$ implies
$$ r'_v\le d_Y(x'_v, p'_v)\le K_1 r'_v. $$  

If $y\in U_v$, then $d_X(y, p_v)< Kr_v$ and so 
$$d_Y(y', p'_v) < d_Y(x_v', p'_v) \eta(K)\le  K_2r'_v .$$
This and the definition of $r'_v$ show
\begin{equation}
\label{qsa1}
 B_Y(p'_v, r'_v)\sub f(B_X(p_v, r_v)) \sub f(U_v) =U'_v \sub 
f(B_X(p_v, Kr_v)) \sub B_Y(p'_v, K_2r'_v). 
\end{equation}
If $u\sim v$, then $U_u\cap U_v \ne \emptyset$ and $r_u\le Kr_v$.
In particular,
$d_X(p_u, p_v) \le K(r_u+r_v) \le K_3 r_v$ and $d_X(x_u, p_v) \le 
d_X(x_u, p_u)+d_X(p_u, p_v)\le r_u+K_3 r_v \le K_4 r_v.$ Hence 
\begin{equation}\label{qsa2}
r'_u \le  d_Y(x'_u, p'_u) \le d_Y(p'_u, p'_v) + d_Y(x'_u, p'_v)
\le  d_X(x'_v, p'_v) (\eta(K_3) +\eta(K_4)) \le K_5 r'_v.
\end{equation}

Suppose $z\in U_v$. Since $d_Y(x'_v, p'_v)\ge r'_v$,  there exists 
$y\in \{p_v, x_v\}$ such that $d_Y(y',z')\ge r'_v/2$. 
Then $d_Y(y,z)\le 2Kr_v$. If now $x\in X$ is an arbitrary point 
with $d_X(x,z)\ge r_v/K$, then
$$
r'_v \le  2 d_Y(y', z') \le 2 d_Y(x', z') \eta(2K^2) \le K_6 d_Y(x', z').
$$
This implies that
\begin{equation}\label{qsa4}
B_Y(z',  r'_v/K_6) \sub f(B_X(z, r_v/K)) \sub f({\cal A}\hbox{-}\St_K(v))
={\cal A'}\hbox{-}\St_K(v) \for z\in U_v.
\end{equation}
The assertion now follows from the fact that ${\cal A}$ is a $K$-approximation
and (\ref{qsa1})--(\ref{qsa4}).  \qed

\begin{lemma}\label{weaklyup}
Suppose $(Z,d)$ is a connected metric space and $((V,\sim),
p, r, {\cal U})$ is a
$K$-approximation of $Z$.  Suppose $L\ge K$ and 
 $W\sub V$ is a maximal  set of combinatorially $L$-separated  vertices.
Then $M=p(W)\sub Z$ is   weakly $\la$-uniformly perfect with 
$\la$ depending only on $L$ and $K$. 
\end{lemma}
\proof  
Note that 
$$ K^{-k(u,v)} \le \frac {r(u)}{r(v)} \le K^{k(u,v)} \for u,v\in V.$$ 
Since $d(p(u),p(v)) \le K(r(u)+r(v))$ whenever $u,v\in V$ with 
$u\sim v$, we obtain 
$$ d(p(u),p(v)) \le 2 r(u) k(u,v) K^{1+k(u,v)} \for u,v\in V.
$$ 
Let  $\la=16L^2K^{4+2L}$. Suppose $w_0, w_1\in W$  such 
 that for $z_0=p(w_0)$
and  $z_1=p(w_1)$ we have that 
$z_0\ne z_1$ and $z_1\in \bar B(z_0, r/\la)$, where 
$0<r\le \diam(M)\le \diam(Z)$.
We claim that  $B(z_0, r)\setminus \bar B(z_0, r/\la)$ contains 
a point in $M$.
Since $w_0\ne w_1$ we have $k(w_0, w_1) \ge L\ge K$ and so 
$U_{w_0}\cap U_{w_1}=\emptyset$ by property (3) of a $K$-approximation.
This implies 
\begin{equation}\label {wup0}
r(w_0)\le  d(z_0,z_1) \le r/\la. 
\end{equation} 
Since $\la >4$ there exist points in $Z$ outside
$B(z_0, r/\sqrt\la)$. The connectedness of $Z$ then implies that there
actually exists $z\in Z$ with $d(z_0, z)= r/\sqrt\la$.
Since ${\cal U}$ is a cover of $Z$, we have
$z\in U_v$ for some $v\in V$.  
Then 
\begin{equation} \label{wup1} 
r(v)\le Kr/\sqrt\la. 
\end{equation} 
For otherwise,
$$\dist(z_0, U_v)\le d(z_0,z)=  r/\sqrt \la < r(v)/K,$$ and so 
$z_0\in N_{r(v)/K}(U_v)\sub \St_K(v)$. This implies $k(w_0, v)\le 2K$ 
which leads to 
$$ r(w_0)\ge K^{-2K} r(v) \ge K^{1-2K} r/\sqrt \la   > r/\la, $$
 contradicting (\ref{wup0}).       

Since $W$ is a maximal $L$-separated set in $V$, there exists
$w_2\in W$ such that $k(w_2, v)<L$.  Let $z_2=p(w_2)\in M$. 
We claim that $d(z_2, z_0)>  r/\la$. Otherwise,
$ d(z_2, z_0)\le r/\la$. If $w_2\ne w_0$, then similarly as above we conclude  
$ r(w_2) \le r/\la$.  But by (\ref{wup0}) this is also true if 
$w_2=w_0$. 
Hence we get in this case 
\begin{eqnarray*}
r/\sqrt\la=d(z_0, z) &\le & d(z_0, z_2) + d(z_2,p(v)) +d(p(v), z) \\
&\le & r/\la +   r(w_2) 2LK^{L+1} + K r(v) \\
&\le & r/\la + (2LK^{L+1} + K^{L+1}) r(w_2)\\ 
&\le & (1+2LK^{L+1}+K^{L+1})r/\la < r\sqrt \la. 
\end{eqnarray*}
which is a contradiction.  

Moreover, by (\ref{wup1})  
\begin{eqnarray*}
 d(z_0, z_2) &\le&   d(z_0, z) + d(z, p(v)) + d(p(v), z_2)\\
&\le & 
r/\sqrt \la +  K r(v) +    2LK^{L+1}  r(v) \\
&\le &   (1+K^2+2LK^{L+2}) r/\sqrt \la \, <\,  r.  
\end{eqnarray*}
 This shows that the point $z_2\in M$ 
is contained in $B(z_0, r)\setminus \bar B(z_0, r/\la)$. 
\qed 

\subsection{Circle packings}
\label{circlepackingsection}
 
We will consider  graphs $G$
  embedded in a metric space $Z$. In this context
we will  regard  $G$ as  a topological space  by identifying each edge of
$G$ with a  unit interval $I:=[0,1]$ and  gluing these intervals according
to the incidence pattern  of the edges.
An {\em embedding} of $G$ into
$Z$ is then just a map of this topological space
into $Z$ which is a homeomorphism onto its image.
 
If the graph $G$  is embedded in $Z$
we will   identify $G$ with its image under
the embedding.  This image is viewed as  a subset of $Z$
with  certain points and arcs distinguished as vertices and
edges, repectively, so that  their incidence pattern is the same as the
incidence pattern of the graph.

Suppose the  graph $G$ 
  is the $1$-skeleton of
a  triangulation $T$ of a topological $2$-sphere.
By the Andreev-Koebe-Thurston
circle packing theorem (cf.\ for example \cite{marrod})
 the graph $G$ can be  realized
as the incidence graph of a circle packing. This means the following. 
Let $V$ be the vertex set of $G$ with the associated incidence 
relation $\sim$. 
 Then 
 there is a family ${\cal C}$
of  pairwise disjoint open nondegenerate
spherical discs $C_v$, $v\in V$, such that
$\bar C_u\cap \bar C_v \ne \emptyset$ for $u,v\in V$ if and  only
if $u\sim v$.
 
We can  always assume that the circle packing is
{\em normalized}. By this we mean that among the centers of the
discs of the circle packing, there are three {\em normalizing}
points which lie on a great circle of $\sph$ and are equally spaced.
A normalization of a circle packing  can always be achieved
by replacing the original circles by their images under a suitably chosen
M\"obius  transformation. To see this note that for three discs with pairwise
disjoint interior there exists a circle orthogonal to the boundary circles
of the discs. This circle can be mapped to a great circle of the 
sphere by a M\"obius  transformation.
The images of the original discs are discs with centers on 
this great circle.  By applying an additional  M\"obius  transformation
fixing the great circle as a set, we can achieve that the centers of the 
discs are equally spaced.

It is easy to see that in a  normalized circle packing all
 discs are smaller than hemispheres. In particular, if two
different discs in the packing have a  common boundary point, then there is a
unique geodesic joining the centers. If we join the centers of
adjacent discs in the circle packing in this  way, then we get an embedding
of $G$ on the sphere. The closures of the complementary
regions of this embedded graph are closed spherical
triangles $\Delta$  forming a  triangulation
$T'$ of $\sph$ combinatorially equivalent to $T$.
If $v\in V$ let $p(v)$ be the center of the disc $C_v$ 
corresponding to $v$, and let $r(v)$ be the spherical  radius of $C_v$.
Let  $U_v$  be the interior  of the union
 of all triangles $\Delta\in T'$  having $p(v)$ as a vertex.
Then $U_v$ is open, starlike with respect to 
$p(v)$ and contains  $C_v$.
Moreover, the  sets $U_v$, $v\in V$,  
 form a cover ${\cal U}$ of $\sph$.
 
Given these definitions we claim:

\begin{lemma} \label{circpack}
Suppose $G$ is combinatorially equivalent to a 
 $1$-skeleton of a triangulation of  $\sph$, 
 and ${\cal C}$ is a normalized
 circle packing realizing $G$.
Then $(G,p,r, {\cal U})$ is a $K$-approxi\-mation of $\sph$ with $K$
depending only on the  valence of $G$.
\end{lemma}
 
\proof It is a well-known fact that for a circle packing
of Euclidean circles the ratio of the radii of two adjacent discs in the packing is bounded
by a constant  depending only on the number of neighbors of (one of)
these  discs (this is called the ``Ring Lemma").  For a packing of 
spherical circles a similar statement is true if no  disc in the packing
is larger than a hemishere, in particular if the packing is normalized. 
In other words, if $u,v\in V$ and  $u\sim v$,
 then $C^{-1} \le r_u/r_v\le C $
with $C$  depending only on the  valence of $G$.
Choosing $K$ suitably depending on  the valence of $G$,
 it is easy to see that the conditions (1)--(5) of  a $K$-approximation
are true for $(G,p,r, {\cal U})$.
We omit the details. \qed

\subsection{Construction of good graphs}

In the following we will work with a modification  of the $LLC_1$-condition
for a metric space $(Z,d)$:

\smallskip\no
($\la$-$\widetilde{LLC}_1$) $\quad$   If $x, y\in Z$, $x\ne y$,  then
there exists an arc $\ga$ with endpoints $x$ and $y$ such that 
$$ \diam(\ga) \le \la d(x,y). $$   

Obviously, $\la$-$\widetilde{LLC}_1$ implies $(1+2\la)$-$LLC_1$.  
A similar quantitative implication in the other direction
  will not be true in general, unless 
$Z$ is locally ``nice". For example, if $Z$ is locally Euclidean, then
a simple covering argument shows that $\la$-$LLC_1$ implies 
$3\la$-$\widetilde{LLC}_1$. So for topological manifolds $LLC_1$ and 
$\widetilde{LLC}_1$  are quantitatively equivalent.  

\begin{lemma}
\label{goodgraphlem} Suppose $(Z,d)$  is a metric space which is 
$C_0$-doubling and $\la$-$\widetilde{LLC}_1$.
Let  $0< r \le \diam(Z)$ and suppose $A\sub Z$ is a maximal $r$-separated
set.  Then  
 there exists
a connected graph $\Ga=(V,E)$ which is embedded in $Z$ and has the following 
properties: 
\begin{itemize} 
\item[{\rm (i)}] The  valence of $\Ga$ is bounded by $K$. 

\item[{\rm (ii)}] The vertex set $V$  contains $A$. 

\item[{\rm (iii)}] If $u,v\in A$ with $d(u,v)<2r$, then $\Ga$ contains an 
arc $\ga$ joining $u$ and $v$
with $\diam(\ga) \le Kr$. The graph $\Ga$ consists
of the union of these arcs.

\item[{\rm (iv)}] For all balls $B(a,r)\sub Z$ we have 
$\# (B(a,r)\cap V) \le K$. 
\end{itemize} 
Here the constant $K\ge 1$ depends only on $C_0$ and $\la$.  
\end{lemma}

\proof 
 For all two-element 
subsets $\{u,v\} \sub A$ with $d(u,v)<2r$ choose an arc $\alpha$
 with endpoints $u,v$ and $\diam(\alpha) \le 2 \la r$. 
Let ${\cal A}$ be the family of arcs thus obtained.
Since $Z$ is doubling, there exists $N_1=N_1(C_0,\la)\in \N$ such that 
each arc in ${\cal A}$ can be covered by at most $N_1$ open balls
of radius $r$.  

We claim that 
there exists $N=N(C_0,\la)\in \N$ such that ${\cal A}$ can be written as 
a disjoint union ${\cal A}={\cal A}_1 \cup \dots \cup {\cal A}_N$, where 
each of the subfamilies ${\cal A}_i$ has the property that if 
$\alpha,\alpha'\in {\cal A}_i$ are two distinct arcs, then 
\begin{equation}\label{arcsep}
 \dist(\alpha,\alpha') >  8\la r. 
\end{equation}

To see that this can be done,  note first that since $Z$ is $C_0$-doubling
there exists $N_2=N_2(C_0,\la) \in \N$ such that 
$$ \#(\bar B(a,12\la r)\cap A)< N_2\for a\in Z. $$
Hence if $\alpha\in {\cal A}$, then 
\begin{equation}\label{arcbd}
 \#\{\alpha'\in {\cal A} : \dist(\alpha, \alpha')\le 8\la r \}
 < N_2(N_2-1)/2.
\end{equation} 
Let $N=N_2(N_2-1)/2$. An argument using Zorn's lemma and (\ref{arcbd}) shows
that there exists a labeling of the arcs in ${\cal A}$ by the numbers
$1, \dots, N$ such that no  two distinct arcs $\alpha, \alpha'\in 
{\cal A}$ with $\dist(\alpha, \alpha')\le 8\la r$ have the same label.
Define ${\cal A}_i$ to be the set of all arcs with label $i$. 

Now define  graphs $\Gamma_i=(V_i,E_i)$  for  $i=1,\dots, N$ 
inductively as follows.  The graphs $\Gamma_i$ will be embedded in 
$Z$,  their edges will have diameter bounded by $2\la r$ and 
we have 
\begin{equation}\label{edgebd}
M_i:=\max_{a\in Z}\,  \#\{  e\in E_i : e \cap B(a, r) \ne \emptyset
  \} \leq (2N_1+4)^i.  
\end{equation} 

Let $\Gamma_1$ be the union of the arcs in  ${\cal A}_1$, where 
we consider these arcs as the edges of $\Gamma_1$ and the set of their
endpoints as the set of vertices.  Note that by 
(\ref{arcsep}) the graph $\Gamma_1$ is embedded   in $Z$ and by choice of the 
arcs in ${\cal A}$ 
  the diameter of each edge will be bounded by  $2\la r$.   Moreover, each ball
$B(a,r)$ can only meet at most one arc in ${\cal A}_1$, so (\ref{edgebd}) is
true  for $i=1$. 

Suppose $\Gamma_{i-1}$ has been constructed. We consider an arbitrary 
arc $\alpha\in {\cal A}_i$ and modify it  to obtain an arc  with 
the same endpoints such that 
for  each edge  $e\in E_{i-1}$ 
the set $\alpha \cap e$ is connected. Note first that the number of edges 
$\alpha$ meets is  bounded by $N_1M_{i-1}$, and in particular finite.
This follows 
from the  definition of $N_1$ and 
$M_{i-1}$.

Let $e_1,\ldots, e_k\in E_{i-1}$ be  the edges that  meet $\alpha$.
Assume inductively that we have modified $\alpha$ into an arc
(also called $\alpha$ by abuse of notation)  such 
that
\begin{equation}\label{connected}
 \text{the sets}\quad 
\alpha\cap e_1,\, \dots, \, \al\cap e_{j-1} \quad\text{are connected.}
\end{equation} 
Let 
$\gamma$ be the smallest (possibly degenerate) subarc  of $\alpha$ 
which contains $\alpha \cap e_j$. Then the endpoints of $\gamma$ 
are contained in $e_j$,  and $\alpha\setminus \gamma$ is disjoint
from $e_j$. Replace $\gamma\sub \alpha$ by the subarc of $e_j$ which has 
the same endpoints as $\gamma$. This new curve $\alpha$ is an arc and the set 
$\alpha\cap e_j$ is connected.  Since the edges in $E_{i-1}$ are nonoverlapping
(i.e., they have disjoint interiors), the statement (\ref{connected}) is still
true for the new arc $\alpha$ (some of the intersections in (\ref{connected})
may have become empty) and there are no new edges that $\alpha$ meets. 
After at most $k$ modifications, the arc $\alpha$ will have the same 
endpoints as before,
and its diameter will be bounded 
by $2\la r+ \sup_{e\in \Ga_{i-1}} \diam(e)\le 4\la r$. 
The arc $\alpha$ has a subdivision into nonoverlapping subarcs
which consists of the sets $\alpha\cap e$ for $e\in E_{i-1}$
 and its complementary subarcs.
Hence $\alpha$ is subdivided into at most $2k+1\le 2N_1M_{i-1}+1$ 
subarcs which all have diameter bounded by $2\la r$. 
Let $\tilde{\cal A}_i$ be the set of these new arcs $\alpha$.
Then for any two distinct arcs in $\tilde{\cal A}_i$ we have  
\begin{equation}\label{arcsep2} 
\dist(\alpha,\alpha') > 2 \la r.  
\end{equation}

The graph $\Ga_i=(V_i,E_i)$ is now obtained   from $\Ga_{i-1}$ and 
the set of modified arcs $\tilde {\cal A}_i$ as follows. If for 
$e\in E_{i-1}$ there exists $\alpha\in \tilde {\cal A}_i$
which meets $e$,  subdivide  $e$ by introducing
new vertices into at most 
three new edges such that $e\cap \alpha$ becomes a vertex or an edge.   
Every edge $e\in E_{i-1}$ is subdivided at most once, since it cannot meet
two distinct arcs in $\tilde {\cal A}_i$ by (\ref{arcsep2}).  
To this graph obtained by subdividing some of the edges of $\Ga_{i-1}$,
we add the edges and vertices from the subdivision of the arcs $\alpha
\in \tilde {\cal A}_i$. Obviously, $\Ga_i$ is embedded in $Z$ and all its edges 
have diameter bounded by $2\la r$. If $B(a,r)$ is an arbitrary 
ball,  then an   edge $e\in E_i$ meeting $B(a,r)$ is either 
a subset of an edge in $E_{i-1}$ meeting $B(a,r)$ or it 
is   an edge  obtained from the subdivision of some 
arc $\alpha\in \tilde {\cal A}_i$.
By (\ref{arcsep2}) all these latter  edges  lie on the same arc $\alpha$.
 Hence 
$$ M_i\le 3 M_{i-1}+ 2N_1M_{i-1}+1 \le (2N_1+4)^i.$$   

We let $\Gamma=\Gamma_N$. Then the underlying set of $\Gamma$ is equal
to the union of the arcs in ${\cal A}_1\cup\tilde {\cal A}_2\cup\dots\cup
\tilde {\cal A}_N$. This shows (ii) and (iii).  These conditions 
imply that $\Ga$ is connected.
Suppose $v$ is a vertex of 
$\Ga$. If an  edge $e$  has a vertex  $v$ as an endpoint, then 
$e\cap B(v,r)\ne \emptyset$.
{}From (\ref{edgebd}) it follows that the number of edges with endpoint
$v$ is bounded by $M_N$ which gives (i). Finally, (iv) follows from 
(\ref{edgebd}) and 
$$ \#(B(a,r) \cap V)\le 2 \#\{e\in E: e\cap B(a,r)\ne \emptyset\}.$$ 
\qed

\begin{proposition}\label{goodtriang0} 
Suppose $(Z,d)$ is a metric space homeomorphic to a $2$-sphere.
If $(Z,d)$ is $C_0$-doubling and $\la$-$LLC$, then for given 
$0<r\le \diam(Z)$  and any maximal $r$-separated set $A\sub Z$ there 
exists an embedded graph $G =(V,E)$
which is the $1$-skeleton of a triangulation
$T$  of
$Z$ such that
\begin{itemize}
\item[{\rm (i)}] The maximal valence of $G$ is bounded by $K$.
 
\item[{\rm (ii)}] The vertex set $V$ of $G$ contains $A$. 

\item[{\rm (iii)}] If $e\in E$, then $\diam(e) < Kr$. If $u,v\in
V$ and $d(u,v)<2r$, then $k(u,v)<K$.    
 
\item[{\rm (iv)}] For all balls $B(a,r)\sub Z$ we have
$\# (B(a,r)\cap V) \le K$.
\end{itemize}
Here the constant $K\ge 1$  depends only on  $C_0$ and $\la$. 
\end{proposition}

Note that since  $G$ is embedded in $Z$, and we can consider the vertices 
and edges of $G$ as subsets of $Z$.
 For $v\in V$ let $p(v)=v$, $r(v)=r$ and $U_v=B(v,Kr)$. Moreover, if 
${\cal U}=\{U_v: v\in V\}$  then under the above assumption we immediately
have: 

\begin{corollary} \label{goodtriang} 
$(G,p,r, {\cal U})$ is a $K'$-approximation of $Z$, where $K'$  depends 
only on $\la$ and $C_0$. 
\end{corollary}

\begin{corollary} \label{goodapprox}
Suppose $Z$ is a metric space homeomorphic to $\sph$. If $Z$ is 
$C_0$-doubling and $\la$-$LLC$, then there exist $K\ge 1$ only depending on
the $C_0$ and $\la$ and a sequence 
${\cal A}_k=(G_k,p_k,r_r, {\cal U}_k)$ of  $K$-approximations  of $Z$, whose 
  graphs $G_k=(V_k, E_k)$ are 1-skeletons of triangulations $T_k$ of $Z$ and
for which 
$$ \lim_{k\to \infty} \mesh({\cal A}_k) =0. $$  
\end{corollary}

\proof This follows immediately from Corollary \ref{goodtriang} if we apply
Proposition \ref{goodtriang0} for a maximal $(1/k)$-separated set $A_k$. 
\qed 

\medskip
\no{\em Proof of Proposition \ref{goodtriang0}.}
First we claim that
 every (continuous)  loop $\phi\: \mathbb{S}^1\ra Z$ 
such that  $\phi(\mathbb{S}^1) \sub B(p,R)$ for some $p\in Z$ and $R>0$ 
 is null-homotopic in $B(p,\la R)$. For this note that   since $Z$ is
$\la$-$LLC$, the compact set $A=Z\setminus B(p,\la R)$ is contained in a
component  of $Z\setminus \phi(\mathbb{S}^1)$. Since $Z$ is homeomorphic
to $\sph$ it follows that $\phi$ is null-homotopic in $Z\setminus A=
B(p,\la R)$.     

Since $Z$ is a topological manifold  and $\la$-$LLC$, it is
$\la'$-$\widetilde{LLC}$ with $\la'=3\la$. 
Let $\Ga_1 =(V_1,E_1)$ be a graph embedded in $Z$ that 
satisfies the conditions (i)--(iv) of Lemma \ref{goodgraphlem} 
 with some constant $K'$ depending
on the data of $Z$.
The idea for constructing $G$ is to subdivide the components of $Z\setminus
\Ga_1$ into triangles. In order that this results in a graph as desired, 
we have to bound the diameter of such a component. We need two lemmas.  

\begin{lemma}
Given a continous map   $f_0\:\mathbb{S}^1\ra Z$, there is a continuous map
$f_1\:\mathbb{S}^1\ra \Ga_1\subset Z$ and a homotopy $f_0\sim f_1$
so that 
 the tracks of the
homotopy have diameter bounded by  $C_1r$ where $C_1$
depends only on $C_0$ and $\la$.
\end{lemma}
\proof
Since $A\sub V$ is a maximal $r$-separated set,  we have 
$\dist(z,A)<r$ for all $z\in Z$. 
Since $f_0(\mathbb{S}^1)$ is compact, for some $r'\in (0,r)$ 
we  have $\dist(f_0(\zeta), A) < r'$ for all $\zeta\in \mathbb{S}^1$.  
 Since $f_0$
is uniformly continuous,   we can find a finite set $S\sub \mathbb{S}^1$   
containing at least two points 
 such that  if $J\sub \mathbb{S}^1-S$ is a maximal 
complementary arc, then $\diam(f_0(J))< r-r'$. 
For  
 each  $\zeta\in S$ we can find a point $f_1(\zeta)\in A$ such that 
$d(f_0(\zeta),f_1(\zeta))< r'$. 
Let $J\sub \mathbb{S}^1-S$ be  a maximal
complementary arc and suppose its endpoints are 
 $\zeta,\zeta'\in S$. 
Then $\dist(f_1(\zeta),f_1(\zeta'))<2r$ and so by 
property (iii) of $\Ga_1$ we can extend $f_1$ continuously to 
$\bar J$ such that $f_1(\bar J)$ is an arc in $\Ga_1$ of diameter at most 
$K'r$.
If we extend $f_1$ in this way to all such arcs $\bar J$, then 
we get a continuous map $f_1\: \mathbb{S}^1\ra \Ga_1$.  

We build a homotopy $H\:\mathbb{S}^1\times I\ra Z$ (where $I=[0,1]$)
  from $f_0$ to $f_1$ as follows.  We set $H(\zeta,0)=f_0(\zeta)$
and $H(\zeta,1)=f_1(\zeta)$ for all $\zeta\in \mathbb{S}^1$.
For each $\zeta\in S$, define $H\restr_{\{\zeta\}\times I}$
to be a  path   connecting 
$f_0(\zeta)$ to $f_1(\zeta)$ of diameter bounded by $\la' r=3\la r$.
 We have defined $H$ on 
$(\mathbb{S}^1\times\{0,1\})\cup(S\times I)$.  If
$J\subset \mathbb{S}^1-S$ is a maximal  complementary arc,
then we can  extend $H$ to $\bar J\times I$ so that
the  image of this set is contained a ball of radius $Cr$ where 
$C=C(C_0,\la)$. Here we use the fact
that the boundary of the ``square'' 
$\bar J\times I$ is mapped into a ball of radius
$R=(3\la +K'+1)r$ and this loop is null-homotopic in a ball with the same
center and radius $\la R$.  
  It follows that the  tracks $t\mapsto H(\zeta, t)$
of the homotopy have diameter bounded by  $C_1r$ with $C_1=C_1(C_0,\la)$.
\qed

\begin{lemma}
The  diameter of each connected component of $Z-\Ga_1$
is bounded by  $C_2r$ where $C_2$  depends only on $C_0$ and $\la$.
\end{lemma}
\proof  We have to show that if $C_2$ is large enough 
depending on the data, then for every point $p \in Z-\Ga_1$
the set $\Ga_1$ separates $p$ and the points in $Z-\Ga_1$ outside 
$B(p, C_2r)$.
Indeed, with the notation of the last lemma
 we can choose  $C_2=4+2C_1$. To see that $C_2$ has the desired 
property    
 we may assume $M=Z\setminus B(p, C_2r)\ne
\emptyset.$  Obviously, $A=\bar B(p,\frac12 (C_2+1)r)\setminus 
B(p, \frac12 (C_2-1)r)$ separates
$p$ from $M$. Using the fact that $Z$ is homeomorphic to a $2$-sphere, it 
is easy to see that there is a Jordan curve 
in an arbitrarily small neighborhood of $A$ separating $p$ from $M$. 
 In particular, there exists a loop $f_0\:\mathbb{S}^1\ra  Z$ 
such that 
$f_0(\mathbb{S}^1)\sub
B(p,\frac12 (C_2+2)r)\setminus  \bar B(p, \frac12(C_2-2)r)$ and 
 the winding number of $f_0$ 
with respect to $p$ differs from the winding number of $f_0$
with respect to any point of $M$.
By the previous lemma we can find a loop $f_1\: \mathbb{S}^1 \ra  \Ga_1$
 homotopic to $f_0$ such that the tracks of the homotopy stay inside 
$$ B(p, {\textstyle\frac12}
(C_2+2+2C_1)r) \setminus  \bar B(p, {\textstyle\frac12}
 (C_2-2-2C_1) r)\sub B(p,Cr)\setminus
\{p\}. $$
  In particular,  the winding number of  $f_1$ with respect 
to  $p$ will still be different from the winding number of 
$f_1$ with respect  to any point in $M$. 
Hence $f_1(\mathbb{S}^1)$ also  separates  $p$ from $M$. 
Since $f_1(\mathbb{S}^1)\sub \Ga_1$,  the point $p$ is separated by $\Ga_1$ 
from the points in $Z\setminus \Ga_1$ outside $B(p, C_2r)$.   
\qed

\medskip
Since $\Ga_1$ is connected, a  component  $\Omega$ of  $Z\setminus \Ga_1$ is 
a simply connected region  whose boundary $\partial \Omega$ is 
a finite union of edges in $\Ga_2$.  Note that by the previous lemma,
the number of these  edges is bounded by a number depending only
on the data of $Z$. 

Now define a new graph $\Ga_2=(V_2,E_2)$ as follows:  Subdivide the edges
of $\Ga_1$ by choosing for each edge  a point in its interior. 
Moreover for each component  $\Omega$ of  $Z\setminus \Ga_1$ choose a point
in its interior. These points together with the the set $V_1$ 
form the vertex set $V_2$ of $\Ga_2$.  The edges 
of $\Ga_2$ are the arcs obtained by the 
subdivision of the edges in $\Ga_1$. Moreover,
for each
component  $\Omega$  of $Z\setminus \Ga$,
we introduce new edges as follows. The vertices in $V_2$ on the boundary 
of $\Omega$ can be brought into a natural  cyclic order $v_1, \dots, v_{N},
v_{N+1}=v_1$, possibly with repetitions, such that successive vertices are
adjacent, i.e., endpoints of an arc obtained from the subdivision of 
the edges in $\Ga_1$.  Note that each vertex can occur at most twice
in this  given cyclic order.  Hence 
 $N$ is bounded by a number depending only on
the data.  Since $\Omega$ is simply connected, we can    connect 
the vertex $v$ chosen in the interior of $\Om$ with each of the vertices 
$v_i$ by an arc $e_i$ such that $e_i\setminus \{v_i\}\sub \Omega$ and such 
that two of these arcs have only the point $v$ in common. 

The graph $\Ga_2$ is embedded in $Z$, and has complementary regions 
whose closures are topological triangles, i.e., there are exactly
 three different  vertices 
and  edges in successive order on the boundary of such a region.
One of these vertices is a vertex contained in  $Z\setminus \Ga_1$, one
will be in the interior of an edge $e\in E_1$ and one vertex will be also
a vertex of $\Ga_1$. 
In particular, the components of $Z\setminus \Ga_2$ are Jordan regions.
In  general, the set of these  triangles which are the closures of 
components of $Z\setminus \Ga_2$ will not be a triangulation of 
$Z$, because it may happen that two such triangles have the same vertex
set without being identical.  This situation arises 
from components of $Z\setminus \Ga_1$ which are not Jordan regions. 

Define a graph $G=(V,E)$ obtained from $\Ga_2$ in the same way as  
 $\Ga_2$ was obtained from $\Ga_1$. Then the closures of the complementary 
components of $Z-G$ are topological triangles which triangulate 
$Z$ so that the $1$-skeleton of this triangulation is $G$. 
The other desired properties of $G$ follow immediately 
from the previous lemma and the properties of $\Ga_1$.  \qed

\subsection{Modulus}
\label{modulus}

Suppose  $(Z,d,\mu)$ is a metric measure  space, i.e., $d$ is a 
metric and $\mu$ a Borel measure on $Z$. Moreover, we assume that $\mu$ 
is locally finite and has dense support. 
The space $(Z,d,\mu)$ is called ({\em Ahlfors})
$Q$-{\em regular}, $Q > 0$, if the measure $\mu$ satisfies
\begin{equation} \label{regular} 
C^{-1}  R^Q \le \mu(B(a,R)) \le CR^Q
\end{equation}
for each open ball $B(a,R)$ of radius $0 < R \le \diam(Z)$ and for some
constant $C \ge 1$ independent of the ball.  The numbers $Q$ and $C$ 
are called the {\em data}  of $Z$. 
If (\ref{regular}) is true for some measure $\mu$, then a similar
inequality holds for $Q$-dimensional Hausdorff measure ${\cal H}^Q$.
Hence, if in a $Q$-regular space the measure is not specified, 
then we assume that the underlying measure $\mu$ is  the 
Hausdorff measure ${\cal H}^Q$. 

We call a Borel function $\rho\: Z \ra  [0,\infty]$ an {\em upper gradient}
of a function $u\: Z\ra \R$ if 
$$ |u(x)-u(y)| \le \int_\ga \rho\, ds, $$
whenever $x,y\in Z$ and $\ga$ is a rectifiable curve joining $x$ and $y$. 
Here integration is with respect to arclength on $\ga$. 

Suppose $B=B(a,r)$ is an open ball in $Z$. If $\la>0$ we let
$\la B:=B(a,\la r)$. Moreover, if $u\: Z\ra \R$  is a locally 
integrable  function on 
$Z$, we denote by $u_B$ the average of $u$ over $B$, i.e., 
$$ u_B= \frac 1{\mu(B)}\int_B u\, d\mu.$$
The metric measure space is said to satisfy a {\em $(1,Q)$-Poincar\'e
inequality}, where $Q\ge 1$, if there exist  constants $C>0$ and 
$\la \ge 1$ such that 
$$ \frac 1{\mu(B)}\int_B |u-u_B|\, d\mu \le C (\diam(B)) \left(
\frac 1{\mu(\la B)}\int_{\la B}  \rho^Q\, d\mu\right)^{1/Q}, $$
whenever $u$  is a locally 
integrable  function on $Z$, the function  $\rho$ is an upper gradient 
of $u$, and $B$ an open ball in $Z$.

A {\em density} (on $Z$)  is a Borel function $\rho\:Z\ra [0,\infty]$.
A density $\rho$ is called  {\em admissible} for a curve family 
$\Gamma$  in $Z$, if 
$$      \int_\ga \rho\, ds \ge 1                $$
for each locally rectifiable curve
  $\ga \in \Gamma$. Here integration is with respect to arclength 
on $\ga$.  If $Q\ge 1$, 
the $Q$-{\it modulus} of a family
$\Gamma$ of curves in $Z$ is the number
\begin{equation}
\Mod_Q (\Gamma) = \inf \int \rho^Q \, d\mu ,
\end{equation}
where the infimum is taken over all densities  $\rho \: Z \ra [0,\infty]$
that are admissible for $\Gamma$. 
If $E$ and $F$ are  (nondegenerate) continua in $Z$, 
we let $\Mod_Q(E,F)$ denote the $Q$-modulus of 
 the family of  curves in $Z$ connecting 
$E$ and $F$. 

Suppose $Z$ is a rectifiably connected metric measure space. 
Then $Z$ is called a $Q$-{\it Loewner space}, $Q \ge  1$,  if 
there exists a positive decreasing  function $\Psi\: \R^+\ra \R^+$ such that 
\begin{equation} 
\Mod_Q (E,F) \ge \Psi(\Delta(E,F))
\end{equation} 
whenever $E$ and $F$ are 
disjoint
continua in $Z$.  Recall  that  $\Delta(E,F)$ is the relative distance of $E$
and $F$ as defined in (\ref{reldist}). 
The number $Q$ and the  function $\Psi$ are the data of the Loewner space
$Z$.  

The Loewner condition was introduced in  \cite{heinkosk} and 
 quantifies the idea that a space has many rectifiable
curves. According to Thm.\ 5.7 and Thm.\ 5.12 in \cite{heinkosk}
a proper $Q$-regular metric space $Z$
satisfies a $(1,Q)$-Poincar\'e inequality if and only if $Z$ is $Q$-Loewner
(Note that the assumption of $\varphi$-convexity in \cite[Thm.\ 5.7]{heinkosk}
is unnecessary, since a proper $Q$-regular metric
space satisfying a   $(1,Q)$-Poincar\'e inequality is quasiconvex
\cite[Appendix]{CHEEGER}).

We will use the following fact about Loewner spaces.

\begin{proposition}
\label{heinkosk1}
Suppose $(Z,d,\mu)$ is a $Q$-regular $Q$-Loewner space, $Q>1$.
Then there exist  constants $\la\ge 1 $ and $C>0$  depending only
on the data of $Z$ with the following property. 

If $z\in Z$, $0<s\le \diam(Z)/\la$,
 and  $Y_1,Y_2\subset Z$
are continua with  $Y_i\cap B(z,s)\neq\emptyset$ and
$\diam(Y_i)\ge s/4$ for $i\in \{1,2\}$,
then for every Borel function $\rho:Z\ra [0,\infty]$
there exists 
a rectifiable curve $\eta$ in  $Z$ joining
$Y_1$ and $Y_2$ such  that
 
$$\int_\eta \rho\, ds\leq
C\left( \int_{B(z,\la s)}\rho^Q\, d\mu\right)^{1/Q}.$$
\end{proposition}
We will skip the proof of this proposition which is very similar to 
the proof of Lem.\ 3.17 in \cite{heinkosk}. Essentially the result is true,
because the relative distance of $Y_1$ and $Y_2$ is bounded by a fixed 
constant. Hence the  regularity and the Loewner condition imply that if 
$\la$ is large enough depending on the data, then the modulus of the 
family of curves   inside $B(z,\la s)$ joining $Y_1$ and $Y_2$ is 
  bigger than a constant.

Suppose $G=(V,\sim)$ is a graph, and
$A,B$ are  subsets of $V$.
We will define the {\em combinatorial $Q$-modulus}  $\GMod^G_Q(A,B)$
of the pair 
$A$ and $B$ as follows. Call a weight function $w\: V \ra [0,\infty]$
{\em admissible}
for the pair $A$ and $B$, if 
$$
\sum_{i=1}^n w(x_i) \ge 1, $$
whenever $x_1, \dots, x_n$ is chain connecting 
$A$ and $B$. 

Now let 
$$ 
\GMod^G_Q(A,B) = \inf \sum_{v\in V} w(v)^Q,
$$
where the infimum is taken over all weights $w$ that are admissible for 
$A$ and $B$. Note that $\GMod^G_Q(A,B) \ge 1$ if $A\cap B\ne \emptyset$. 
We drop the superscript $G$ in $\GMod^G_Q(A,B)$ if the graph $G$ is understood.

If $A\sub V$ and $s> 0$ we  denote by 
$N_s(A)$ the $s$-neighborhood of $A$, i.e., the set of all $u\in V$ 
for which there exists $a\in A$ with $k_G(a,v)< s$. 

If we want to estimate the  $Q$-modulus of the pair $(A,B)$, then the following
lemma will allow us to change   the sets $A$ and $B$ with quantitative 
control.

\begin{lemma} \label{setchange} 
Suppose $G=(V,\sim)$ is a graph with valence  bounded by $d_0\ge 1$.
For every $Q>1$ and $s> 0$ there exists a number $C=C(d_0, s, Q)$ with 
the following property: If $A,B, A', B'\sub V$, $A'\sub N_s(A)$, and 
$B'\sub N_s(B)$, then 
$$ \GMod_Q(A',B')\le C \GMod_Q (A,B).  $$ 
\end{lemma} 

\proof Note that if $w$ is admissible for $A$ and $B$, then 
$\tilde w\: V \ra [0,\infty]$ defined by 
$$ \tilde w(v) =\sum_{u\in B(v,s)} w(u) \for v\in V$$ 
is admissible for  $(A',B')$. Moreover, since the vertex degree of 
$G$ is bounded, it follows that each ball $B(v,s)$ has a cardinality
bounded by a constant  depending only on $s$ and $d_0$. 
It follows that 
$$\sum_{v\in V} \tilde w(v)^Q\le C \sum_{v\in V} w(v)^Q,$$
with $C=C(s,d_0, Q)$. 
The lemma follows. \qed

\subsection{$K$-approximations and modulus comparison}
 
In this section we
relate the $Q$-modulus on  a metric  space 
to the $Q$-modulus on the graph of a  $K$-approximation.
Results of this
general nature  are well-known. 
The (minor) novelty here is that the local scales  may vary from 
point to point.

Let $(Z,d)$ be a metric space. Throughout this section 
${\cal A}=(G,p,r,{\cal U})$ will be a  $K$-approximation of $Z$ with 
graph $G=(V,\sim)$.
For each subset $E\subset Z$ we define $V_E\defeq\{v\in V\mid
U_v\cap Z\neq\emptyset\}$. Note that $V_E\sub V$ depends on  ${\cal A}$, 
but we suppress this dependence in our notation. 

\begin{proposition}\label{xferprop}
Let $(Z,d,\mu)$ be  a $Q$-regular metric measure  space, $Q\ge 1$, and
let ${\cal A}$ be a $K$-approximation of $Z$.
Then there exists a constant $C\ge 1$ depending only on $K$ and the data
of $Z$ with the following property:
 
 If $E,F\sub Z$ are continua 
 and if $\dist(V_E, V_F) \ge 4K$, then
\begin{equation} \label{trans1}
\Mod_Q(E,F) \le C
\GMod_Q (V_E,V_F).
\end{equation}
\end{proposition}

\proof
Let $w\:V\ra [0,\infty]$ be an admissible function
for the pair $(V_E,V_F)$: if $v_1\sim\dots\sim v_k$ is a chain 
in $V$ with $v_1\in V_E$ and $v_k\in V_F$, then
$\sum_{i=1}^k w(v_i)\geq 1$.
Define $\tilde w\:V\ra\R^+_0$ by the formula 

$$\tilde w(v)=\sum_{u\in B(v,K)}w(u),$$
and 

$$\rho\defeq \sum_{v\in V}\left(\frac{\tilde w(v)}{r_v}\right)
\chi_{{}_{\scriptstyle{\St_K(v)}}}$$
where $\chi_{{}_{\scriptstyle{Y}}}$ denotes the characteristic 
function of $Y\subset Z$.

\medskip
\no
{\em Mass bounds for $\rho$.}  The cover 
$\{\St_K(v): v\in V\}$ has controlled overlap  
depending on $K$ and there exists a constant $C=C(K)$ such that
$\St_K(v)\sub B(v, Cr_v)$ for $v\in V$. Moreover,  $Z$ is
$Q$-regular and 
every $K$-ball in $V$ has
cardinality controlled by $C(K)$. So  we have
that 

\begin{equation} \label{mass1}
\begin{array}{lcl}\displaystyle \int_Z\rho^Q\, d\mu &\lesssim & 
\displaystyle
\sum_{v\in V}\int_Z\left
(\frac{\tilde w(v)}{r_v}\chi_{{}_{\scriptstyle{\St_K(v)}}}
\right )^Q\, d\mu  \\ &&\\
&\lesssim  & \displaystyle
\sum_{v\in V}\tilde w(v)^Q \lesssim  \sum_{v\in V}w(v)^Q. 
\end{array}
\end{equation}

\medskip
\no
{\em Admissibility of $\rho$.} Now let $\ga\:J\ra Z$  be a rectifiable
path connecting $E$ to $F$. Since ${\cal U}$ is a cover 
of the path  $\ga$, there exists a set $W=\{v_1, \dots, v_k\}$ 
in $V$ such that $\ga \cap U_{v_i}\ne \emptyset $ for 
$i\in \{1,\dots, k\}$,  
 $U_{v_{i}}\cap U_{v_{i+1}}\ne \emptyset$ for $i\in \{1,\dots, k-1\}$,  
and $v_1\in V_E$ and $v_k\in V_F$. The combinatorial distance 
of $v_i$ and $v_{i+1}$ is less than  $K$. Hence there exists a  chain 
$A$ in $V$ connecting $V_E$ and $V_F$ satisfying
$W\sub A\sub N_K(W)$.     

For each $v\in W$, let $J_v:=\ga^{-1}(\St_K(v))$ and 
 $\ga_v:= \ga\restr_{J_v}$.
 Then
the definition of $\rho$ gives
$$\rho(\ga(t))\ge \tilde w(v)/ r_v\for t\in J_v. $$
By our assumption that $\dist(V_E,V_F)\ge 4K$ the path  $\ga$ is not contained 
in any $K$-star of a vertex. For if $\ga\sub \St_K(u)$, 
then there exist $u_1, u_2\in V$ with $k(u_1, u)< K$, $k(u_2, u)< K$, 
$U_{v_1} \cap U_{u_1}\ne \emptyset$, and 
$U_{v_k} \cap U_{u_1}\ne \emptyset$. Then 
$k(v_1, u_1) < K$ and $k(v_k, u_2) < K$ which implies
$\dist(V_E,V_F)\le k(v_1,v_k)<4K$.

Since  $\ga$ is not contained
in any $K$-star of a vertex, we have that
  if a set $U_v$ meets $\ga$, then 
$\length(\ga\cap \St_K(v))\geq r_v/K$
 by  condition (4) of a $K$-approximation. In particular,
for each $v\in W$ we have 
$\length(\ga_v)\geq r_v/K$, and  so

$$\int_{\ga_v}\rho \, ds \gtrsim 
\left (\frac{\tilde w(v)}{r_v}\right )\length(\ga_v)
\gtrsim \tilde w(v).$$
Hence

$$\sum_{v\in W}\int_{\ga_v}\rho\, ds
\gtrsim \sum_{v\in W}\tilde w(v)
\gtrsim \sum_{v\in N_K(W)}w(v) \ge 1, 
$$
since $N_K(W)$ contains the chain $A$ connecting $V_E$ and $V_F$
 and $w$ 
is admissible.  The sets $\St_K(v)$ and hence the sets
    $J_v \sub J $ for $v\in W$
have 
 controlled overlap  depending on $K$, giving

\begin{equation}
\label{admissible1}
\int_\ga\rho\, ds\gtrsim \sum_{v\in W}\int_{\ga_v}
\rho\, ds \gtrsim 1. 
\end{equation}
Combining
(\ref{mass1}) with (\ref{admissible1}) we get
$$\Mod_Q(E,F)\lesssim \GMod_Q(V_E,V_F).$$
\qed

\medskip
It is an interesting question when an inequality like (\ref{trans1}) holds
in the opposite direction. We will not need such a result for the proof
of our theorems, but we will nevertheless explore this question, 
because it illuminates the general picture.
In order to get the desired  inequality, we have to  add
an analytic assumption on  $Z$ to our hypotheses. It suffices to assume 
that $Z$ is a  $Q$-regular $Q$-Loewner space, but as the next proposition 
shows it is enough that a  Loewner type condition holds locally on a scale
 corresponding
to the scale of our $K$-approximation ${\cal A}$.

\begin{proposition}\label{xferprop2}
Let $(Z,d,\mu)$ be  a $Q$-regular metric measure  space, $Q\ge1$, and
let ${\cal A}$ be a $K$-approximation of $Z$.

Suppose that there exists constants $c_1,C_1>0$ with the following 
property: Let $v\in V$, $z\in U_v$, and $0<s\le c_1 r_v$. 
If  $Y_1, Y_2\sub Z$  are closed connected   with 
$Y_i\cap B(z, s) \ne \emptyset $ and $\diam(Y_i)\ge s/4$ for 
$i\in \{1,2\}$, then  for every Borel function $\rho\:Z\ra 
[0,\infty]$ there exists a 
rectifiable path $\eta$ connecting $Y_1$ and $Y_2$  such that 
\begin{equation} \label{acond}  
\int_\eta \rho\, ds\leq
C_1\left( \int_{\St_K(v)}\rho^Q\, d\mu\right)^{1/Q}. 
\end{equation}

Then there exists a constant $C\ge 1$ depending only on $K$, the data
of $Z$, and the constants associated to the analytic condition 
(\ref{acond}) with the following property:
 
 If $E,F\sub Z$ are continua not contained in any $K$-star $\St_K(v)$,
$v\in V$, then
\begin{equation} \label{trans2}
\GMod_Q (V_E,V_F) \le C \Mod_Q(E,F).
\end{equation}
\end{proposition}

Note that by Proposition \ref{heinkosk1} and by the properties
of a $K$-approximation 
every $Q$-regular $Q$-Loewner space $Z$ 
with $Q>1$ satisfies the analytic condition (\ref{acond}) with appropriate
constants  depending only on $K$ and the data of $Z$. 
So Proposition \ref{xferprop} and Proposition \ref{xferprop2} together imply 
the following corollary.  

\begin{corollary}\label{xferprop3}
Let $Z$ be  a $Q$-regular $Q$-Loewner  space, $Q>1$, and
let ${\cal A}$ be a $K$-approximation of $Z$.
Then there exists a constant $C\ge 1$ depending only on $K$ and the data
of $Z$ with the following property:

If $E,F\sub Z$ are continua not contained in any $K$-star 
 and if $\dist(V_E, V_F) \ge 4K$, then
\begin{equation} \label{trans3}
C^{-1} \Mod_Q(E,F)\le \GMod_Q (V_E,V_F) \le C \Mod_Q(E,F).
\end{equation}
\end{corollary}

\bigskip
\no
{\em Proof of Proposition \ref{xferprop2}.\ }
Let $\rho\:Z\ra\R_+$ be an admissible Borel function
for the pair $(E,F)$, i.e.\ 

$$\int_\ga \rho\, ds\geq 1$$
for any rectifiable curve $\ga$ joining $E$ with $F$.
Define $ w\:V\ra\R_+$ by 
$$w(v)\defeq \left (\int_{\St_{3K}(v)} \rho^Q \, d\mu \right )^{1/Q}.$$
\medskip
\no
{\em Mass bound for $w$.}
 Since the numbers  $\#B(v, 3K)$ for $v\in V$ and the 
  multiplicity of the cover  ${\cal U}$
are bounded   by  a constant  depending only on $K$, 
 we have
\begin{eqnarray} 
\sum_{v\in V}w(v)^Q
& \lesssim & \sum_{v\in V} \sum_{u\in  B(v,3K)} \int_{U_u} \rho^Q\,
d\mu\nonumber
\\
&\lesssim & \sum_{v\in V} \int_{U_v} \rho^Q\, d\mu \label{mass12}  \\
&\lesssim & \int_Z\rho^Q\, d\mu. \nonumber 
\end{eqnarray}

\medskip
\no
{\em Admissibility of $w$.}  This step in the 
proof is modelled on arguments from \cite{heinkosk}, and
is based on repeated application  
of our analytic condition.  
We  use  this 
near a single set $U_v$ to prove that under our   assumptions
we have:

\begin{lemma}
\label{leminUv} 
Suppose $v\in V$, and $Y_1,\,Y_2\subset Z$ are closed, connected
sets with $Y_i \cap \St_K(v) \ne  \emptyset $, and $\diam(Y_i)\geq
r_v/(2K^2)$.  Then there is a rectifiable curve
$\eta$ connecting $Y_1$ and $Y_2$  such that 
\begin{equation}
\label{inUv}
\int_\eta \rho\, ds\leq Cw(v), 
\end{equation}
where $C>0$  depends only on $K$ and the data of $Z$. 
\end{lemma}
\proof
Pick $z_1,\,z_2\in \St_K(v)$ so that $z_i\in Y_i\cap \St_K(v)$.
Applying condition (5) of a $K$-approximation repeatedly 
 we find  a curve $\ga$ joining
$z_1$ to $z_2$ so that $\ga\sub \St_{2K}(v)$.
Let 
$$s: = (c_1/K^2) \min_{u\in B(v,2K)} r(u)\simeq r(v), $$   
where $c_1$ is the constant in the hypothesis of Proposition \ref{xferprop2}.
Since $Z$ is $Q$-regular, it is doubling. Moreover, 
$s\simeq r(v)$ and $\diam(\ga)\lesssim r(v)$. 
Hence  the cardinality 
of a maximal $(s/2)$-separated set on $\ga$ is bounded 
by a number  depending only on the data. 
Since $\ga$ is connected, we can find an appropriate subset 
$x_1,\ldots,x_N$ 
of such a maximal set such 
that $d(z_1, x_1)< s$, $d(z_2, x_N)<s$,  and  $d(x_{i-1}, x_i) < s$ for 
$i\in \{2, \dots, N\}$, where $N\in \N$ is bounded 
by a number  depending only on the data.

Now let $\la_1\defeq Y_1$ and $\la_{N+1} \defeq Y_2$. Then 
$\diam(\la_1)\wedge \diam(\la_{N+1}) \ge s/4$ by our assumptions. 
  If $N\ge 2$, we have 
$\diam(\ga)\ge s/2$ and so in addition 
we can find
 continua $\la_i\sub \ga$  with $x_i\in \la_i \sub B(x_i, s)$
and   $\diam(\la_i)\ge s/4$ for $i\in\{2,\dots,N\}$.  

Now  $x_i  \in \ga \sub \St_{2K}(v) $ and so $x_i\in U_{u_i}$ for some
$u_i\in V$ with $k(u_i,v)\le 2K$. Then by definition of $s$ we have
$s\le c_1r_{u_i}$.
Hence we can inductively find rectifiable 
 curves $\eta_1, \dots, \eta_N$ joining 
$\la_1\cup \eta_1\cup \dots\cup \eta_{i-1}$ and $\la_{i+1}$ such that 
\begin{equation}
\label{applheinkosk}
\int_{\eta_i}\rho\, ds \lesssim 
\left( \int_{\St_K(u_i)}\rho^Q\, d\mu\right)^{1/Q} 
\lesssim \left( \int_{\St_{3K}(v)}\rho^Q\, d\mu\right)^{1/Q} =w(v)^Q. 
\end{equation}
This follows from an application of our analytic 
assumption to the ball
$B(x_i,s)$ and  the pair
$\la_1 \cup \eta_1\cup\dots\cup\eta_{i-1} $ and $\la_{i+1}$. 
Note that $\la_{i+1}$ meets $B(x_i, s)$.
The same is true
for the set $\la_1 \cup \eta_1\cup\dots\cup\eta_{i-1}$, since
it meets   $\la_i$ by induction hypothesis.
The union $\eta_1\cup\ldots\cup\eta_{N}$ contains
a rectifiable curve $\eta$ connecting $Y_1$ and $Y_2$ 
  with 
$$\int_{\eta} \rho\, ds\lesssim  N w(v)\simeq w(v).  $$ 
\qed

\medskip
Now suppose $v_1,\ldots,v_k$
are the vertices of a chain  in $G$ joining
$V_E$ to $V_F$.  Then  $ U_{v_1}\cap E\ne \emptyset$,
$U_{v_k}\cap F \ne \emptyset$, and 
 $U_{v_{i-1}}\cap U_{v_i} \ne \emptyset$ for $i\in \{2,\dots,k\}$.   Set
$\la_1\defeq E$, $\la_{k+1}\defeq F$, and for
$i\in \{2,\dots,k\}$ let $\la_i$ be
a closed connected set with  $\la_i\sub \St_K(v_{i-1}) \cap \St_K(v_{i})$ 
and
$$\diam(\la_i)\ge (r_{v_{i-1}}\wedge r_{v_{i}})/ (2K)\ge 
(r_{v_{i-1}}\vee r_{v_{i}})/ (2K^2). $$ 
These sets exist by condition (4) of a
$K$-approximation and the fact that the complement of any 
$K$-star contains elements of $E$ and $F$ and is thus  nonempty.
Moreover,  the fact that $E$ and $F$ are not contained 
in any $K$-star also shows 
$\diam(\la_1)\geq{r_{v_i}}/{K}$
and $\diam(\la_{k+1})\geq{r_{v_i}}/{K}$.

We can inductively find rectifiable curves $\eta_1,\dots, \eta_k$
with 
$$\int_{\eta_i}\rho\, ds\leq C_1w(v_i)$$
so that $\eta_i$ joins $\la_1\cup\eta_1\cup\ldots\eta_{i-1}$
to $\la_{i+1}$. Here $C_1$  depends only on $K$ and the data of $Z$. 
This follows from an application of 
Lemma \ref{leminUv} with $v=v_i$,
  $Y_1:= \la_1\cup\eta_1\cup\ldots\eta_{i-1}$
and $Y_2=\la_{i+1}$. Note that we have
 $\diam(Y_1) \wedge \diam(Y_2) 
\ge r_{v_i}/(2K^2)$. Moreover,  both sets $Y_1$ and $Y_2$ meet
 $\St_{K}(v_i)$.
This  is true for $Y_1=\la_1\cup\eta_1\cup\ldots\eta_{i-1}$,
 since this set meets 
$\la_i$ by induction hypothesis. 

The union $\eta_1\cup\ldots\cup\eta_k$ will contain
a rectifiable curve $\eta$ joining $E$ to $F$
with 

$$1\leq\int_\eta \rho\, ds\leq C_1\sum_{i=1}^kw(v_i).$$
Therefore $C_1w$ is an admissible
test function for $(V_{E},V_{F})$. Hence by (\ref{mass12})

$$\GMod_Q(V_{E},V_{F})\lesssim \Mod_Q(E,F).$$
This completes the proof of Proposition \ref{xferprop}. \qed 

\subsection{The Ferrand cross-ratio}
\label{theferrandcross-ratio} 

If a map quantitatively distorts the modulus
of curve families, then in some situations it follows that the map is \qm. 
A result of this type is the following 
proposition, which illustrates the importance of the concept of a Loewner 
space.  (Cf.\ Remark 4.25 in  \cite{heinkosk}, where a related result 
is mentioned without proof.)

\begin{proposition}
\label{FCR0}
Let   $X$ and $Y$ be metric spaces,
$f\: X\ra Y$ a homeomorphism, and $Q>1$. 
Suppose  $X$ is a $Q$-regular $Q$-Loewner space, $Y$ is $Q$-regular and 
$LLC$, and that 
there exists a constant $K>0$ such that 
\begin{equation} \label{qc}
\Mod_Q(\Gamma) \le K \Mod_Q(f(\Gamma))
\end{equation} 
for every family  $\Gamma$ of curves in $X$.  

Then $f$ is $\eta$-\qm\ with $\eta$ depending only on $K$ and the 
data of $X$ and $Y$. 
\end{proposition} 

Here $f(\Gamma)$ is the family of all curves $f\circ \ga$ with $\ga \in
\Gamma$.

\proof As a Loewner space, $X$ is $\la$-$LLC$ with $\la$ depending 
on the data of $X$, and in particular 
connected. Moreover, $Y$ is $C_0$-doubling
with $C_0$ depending only on the data of $Y$.
So by Lemma \ref{weakqm} it is enough to show that if $(x_1,x_2,x_3,x_4)$ is
a four-tuple of distinct points in $X$ with $[y_1,y_2,y_3,y_4]$ small, 
where $y_i=f(x_i)$, then $[x_1,x_2,x_3,x_4]$ is small,  quantitatively.

Now if $[y_1,y_2,y_3,y_4]$ is small, then by Lemma \ref{goodcont}
we can find continua
$E',F'\sub Y$ with $y_1,y_3\in E'$, and $y_2,y_4\in F'$ such that
 $\Delta(E',F')$ is large, quantitatively. Let $\Gamma'$ be the family 
of all curves in $Y$ joining $E'$ and $F'$, and let $\Gamma$ be the family 
of all curves in $X$ joining $E:=f^{-1}(E')$ and $F:=f^{-1}(F')$.
Then $\Gamma'=f(\Gamma)$ and so by our hypotheses we have
$$ \Mod_Q(E,F)=\Mod_Q(\Gamma) \le K \Mod_Q(\Gamma')= K\Mod_Q(E',F'). $$
Since  $Y$ is $Q$-regular, we have that
$$\Mod_Q(E',F') \lesssim \frac 1 {(\log(1+\Delta(E',F'))^{Q-1}}.$$ 
This is a standard fact following from the upper mass bound for the Hausdorff
measure in $Y$.  It can be 
 be established similarly as  Proposition \ref{FCR2} below.
Hence if $\Delta(E',F')$ is large,  then 
 $\Mod_Q(E',F')$ and so $\Mod_Q(E,F)$ are
small, quantitatively. But in a Loewner space,
we have $$\Phi(\Delta(E,F)) \le \Mod_Q(E,F),$$ where $\Phi\: \R^+\ra  \R^+$
is a positive
decreasing function.  It follows that $\Delta(E,F)$ is large, quantitatively.
Finally, by Lemma \ref{goodcont}
again, this means that for the points $x_1,x_3\in 
E$ and $x_2,x_4\in F$ we have that $[x_1,x_2,x_3,x_4]$ is small,
quantitatively. 
\qed 

\medskip
We will actually  not use this proposition, but rather  corresponding
discrete versions of this result (the closest discrete 
analog is Proposition \ref{unifqm}). We included 
Proposition \ref{FCR0} to  clarify the basic idea.

The relevant point in the preceding proof was   that 
the  cross-ratio of four points 
can be quantitatively controlled by an appropriate modulus.
  So suppose $X$ is a metric  measure  space and let $(x_1,x_2,x_3,x_4)$ be a
four-tuple of distinct points.  For $Q\ge 1$
  define the {\em Ferrand cross-ratio} of the four points to be
\begin{equation} \label{Ferrand}
[x_1,x_2,x_3,x_4]_Q= \inf\,\Mod_Q(E,F),
\end{equation}
where the infimum is taken over all continua $E,F\sub X$ with
$x_1,x_3\in E$ and $x_2, x_4\in F$.
 
Using Lemma \ref{goodcont}, it is not hard to see that if $X$ is a $Q$-regular
$Q$-Loewner space, then the cross-ratio $[x_1,x_2,x_3,x_4]$ is small
if and only if the Ferrand cross-ratio $[x_1,x_2,x_3,x_4]_Q$ is small.
Moreover, if $X$ is only $LLC$ and $Q$-regular, then at least one of these
implication holds. Namely, if $[x_1,x_2,x_3,x_4]$ is small,
then $[x_1,x_2,x_3,x_4]_Q$
is small. The purpose of this section is to establish 
similar results
 for vertices  in a graph coming from a $K$-approximation.

Assume $Q\ge 1$ is fixed and let  $G=(V,\sim)$ be  a connected graph.
Imitating the definition of the Ferrand cross-ratio
 in a metric measure space 
$Z$, we define the Ferrand cross-ratio 
of a four-tuple $(v_1,v_2,v_3,v_4)$ of distinct points in $V$ by 
$$ [v_1,v_2,v_3,v_4]_Q =\inf\,\GMod_Q(A,B), $$
where the infimum is taken over all chains $A,B\sub V$ with $v_1,v_3\in A$ and
$v_2,v_4\in B$. 

\begin{proposition}\label{FCR1}
Let   $Z$
be    a metric measure   space which is $LLC$, let
 ${\cal A}$
be  a   $K$-approxi\-mation of $Z$, and $Q\ge 1$.  
Suppose
that there exists a number $L>0$ and a function $\Psi\:\R^+\ra (0,\infty]$ with
$\lim_{t\to \infty} \Psi(t)=0$ such that
\begin{equation}\label{FCR11}
 \GMod_Q(V_E,V_F) \le \Psi(\Delta(E,F)),
\end{equation}

whenever $E,F\sub Z$ are continua not  contained in any $L$-star.
 
 Then there exists
a function $\delta\:\R^+\ra \R^+$  depending only on $K$, $L$,
 $Q$, $\Psi$ and the
data of $Z$  with the following property:
 
If  $\eps>0$ and  $(v_1,v_2,v_3,v_4)$ is an   arbitrary
four-tuple of vertices in $G$ such that $k(v_i,v_j) \ge 2(K+L) $ for
$i\ne j$, then
$$
[p(v_1), p(v_2), p(v_3), p(v_4)]<\delta(\eps)\Ra
[v_1,v_2,v_3,v_4]_Q< \eps.
$$
\end{proposition}

We will see below (cf.\ Proposition \ref{FCR2}) that if 
$Z$ is  $LLC$
and $Q'$-regular with $Q'\le Q$, then
 condition (\ref{FCR11}) is satisfied 
with  $L=K$ and some function $\Psi$ only depending on $K$ and the data 
of $Z$ (and not on ${\cal A}$).
 
\proof Let $p_i=p(v_i)$ for $i\in\{1,\dots, 4\}$.
Our assumption on the combinatorial
separation of the vertices $v_i$ and properties (2) and (3) of a ${\cal A}$
imply that the points $p_i$ are distinct. Hence $[p_1,p_2,p_3,p_3]$
is well-defined.

We have to show that if $k(v_i,v_j)\ge 2(K+L)$
for $i\ne j$ and  $[p_1,p_2,p_3,p_4]$ is small,
then $[v_1,v_2,v_3,v_4]_Q$ is small, quantitatively.
If $[p_1, p_2, p_3, p_4]$ is small, then by Lemma \ref{goodcont}
there exist continua $E$ and $F$ with $p_1,p_3\in E$, $p_2,p_4\in F$
and $\Delta(E,F)$ large, quantitatively.
Since $E$ is a continuum, we can find a chain $A\sub N_K(V_E)$ connecting
$v_1, v_3\in V_E$. Similarly, we can find a chain $B\sub N_K(V_F)$
connecting $v_2,v_4\in V_F$. Lemma \ref{setchange}
implies  that there exists a
constant $C=C(K)$ such that
$$\GMod_Q (A,B)\le C \GMod_Q  (V_E, V_F).$$

The set  $E\supset \{p_1, p_3\}$
is not contained in the $L$-star of any $v\in V$. For if $E\sub \St_L(v)$,
then there exit  $u_1, u_2 \in V$ with $k(v,u_1)\le L$, $k(v,u_2)\le L$,
$p_1\in U_{u_1}$, and $p_3\in U_{u_2}$. But then 
$p_1\in  U_{u_1}\cap U_{v_1}$ which implies $k(v_1, u_1)<K$ by 
property (3) of a $K$-approximation. Similarly, $k(v_3, u_2)<K$.
Putting these inequality together we get  $k(v_1, v_3)< 2(K+L)$ which
contradicts our assumption on the combinatorial separation of the 
vertices $v_i$.
In the same way we see that $F$  can neither be contained in an $L$-star.
Now from our assumption we obtain 
$$ [v_1,v_2,v_3,v_4]_Q \le \GMod_Q (A,B) \lesssim \GMod_Q  (V_E, V_F)
\le \Psi(\Delta(E,F)). $$
Since $\Delta(E,F)$ is large and 
 $\Psi(t)\to 0$ as $t\to \infty$, this implies that 
$ [v_1,v_2,v_3,v_4]_Q$ is small, quantitatively. 
\qed

\begin{proposition}\label{FCR4}
Let   $Z$
be    a metric measure   space,
 ${\cal A}$ be
 a   $K$-approxi\-mation of $Z$, and
  $t,Q\ge 1$. Suppose
that there exists a number $M>0$ and  a decreasing positive
function $\Phi\:\R^+\ra \R^+$
such that
\begin{equation}\label{FCR41}
 \Phi(\Delta(E,F))\le \GMod_Q(V_E,V_F),
\end{equation}
whenever $E,F\sub Z$ are continua  with $\dist(V_E, V_F)\ge M.$
 
 Then there exists
a function $\delta\:\R^+\ra \R^+$  depending only on  $K$, $M$, $Q$,
 and $\Phi$
 with the following property:
 
If  $\eps>0$ and  $(v_1,v_2,v_3,v_4)$ is   an   arbitrary
four-tuple of  vertices in $G$ such that $k(v_i,v_j)\ge K$
for $i\ne j$, then  we have:
$$
[v_1,v_2,v_3,v_4]_Q< \delta(\eps)\Ra
[p(v_1), p(v_2), p(v_3), p(v_4)]
< \eps.
$$
\end{proposition}

It follows from Proposition \ref{xferprop} that if $Z$ is a $Q$-regular
$Q$-Loewner space, then
 condition (\ref{FCR41}) is satisfied  with 
$M=4K$ and some function $\Phi$ depending only 
on $K$ and the data of $Z$ (and not on ${\cal A}$). 
 
\proof
Let $p_i=p(v_i)$ for $i\in\{1,\dots, 4\}$. Our assumption on the combinatorial
separation of the vertices $v_i$ 
implies that the points $p_i$ are distinct and  $[p_1,p_2,p_3,p_3]$
is well-defined.
 
If  $[v_1,v_2,v_3,v_4]_Q$ is small, then there exist
chains $A,B$ in $G$ with $v_1,v_3\in A$
and $v_2,v_4\in B$ for which $\GMod_Q(A,B)$ is small, quantitatively.
 
We may assume $\dist(A,B) \ge M+4K$. Otherwise, $A'=N_{M+4K}(A)$ and
$B'=N_{M+4K}(B)$  have nonempty intersection which by Lemma \ref{setchange}
leads to
$$ 1\le \GMod_Q(A',B') \le C(K,M,Q) \GMod_Q(A,B). $$
 
Since $A$ is a chain connecting
$v_1$ and $v_3$,  there are elements $u_i$ in $A$
with  $u_1=v_1\sim \dots \sim u_n=v_3$.
Then $U_{u_i}\cap U_{u_{i+1}} \ne \emptyset$
and we can find a curve $\ga_i\sub \St_K(u_i)\cup \St_K(u_{i+1})$
connecting $p(u_i)$ and $p(u_{i+1})$ for $i\in \{1,\dots, n-1\}$.
The  union $E=\ga_1\cup \dots \cup \ga_{n-1}$ is a continuum
joining $p_1$ and $p_3$ with
$$E\sub \bigcup_{i=1}^n \St_K(u_i).$$
If $u\in V_E$, then $U_u\cap  U_w\ne \emptyset $ for some
$w\in N_K(A)$. Hence $V_E \sub N_{2K}(A)$.
A continuum $F$ in $Z$ connecting $p_2$ and $p_4$  with
$V_F\sub N_{2K}(B)$ can be constructed
 in the same  way.
Then $\dist(V_E,V_F)\ge \dist(A,B)-4K\ge M$ and so
$$  \Phi(\Delta(E,F)) \le \GMod_Q(V_E, V_F) \lesssim \GMod_Q(A,B).$$
Since $\GMod_Q(A,B)$ is small, we see that $\Delta(E,F)$ is large,
quantitatively. Lemma \ref{goodcont} implies that $[p_1,p_2,p_3,p_4]$
is small, quantitatively.  \qed

Now we can prove a discrete version of Proposition \ref{FCR0}. 

\begin{proposition} \label{unifqm}
Let $Q\ge 1$, and let   $X$ and $Y$  be    metric  measure  spaces with
$K$-approximations ${\cal A}=(G, p, r, {\cal U})$ and
${\cal A'}=(G, p', r', {\cal U}')$, respectively, whose underlying
graph $G=(V,\sim)$ is the same.
Suppose $X$ is connected,  and $X$ and
 ${\cal A}$ satisfy condition  (\ref{FCR41}) for some $M>0$ and some function
$\Phi$. Suppose $Y$ is $LLC$ and doubling,
  and  $Y$ and ${\cal A}'$ satisfy 
condition (\ref{FCR11})  for some $L>0$ and some  function $\Psi$.
Assume $W\sub V$ is a maximal  set  of vertices with mutual combinatorial
 distance   at least
$s$, where $s\ge  2(K+L)$.  Let $A=p(W)$, $B=p'(W)$ and define
$$ f\: A\ra B, \ x \mapsto p'({p}^{-1}(x)).$$
 
Then $f$ is $\eta$-\qm\ with $\eta$  depending only on $K$, $Q$, $L$, $M$, $s$,
$\Phi$, $\Psi$, and the data of $Y$ (i.e., the parameters in the
$LLC$ and doubling conditions).
\end{proposition}

Since the concept of modulus on a graph is independent of the 
concept of a $K$-approximation, 
the analog of (\ref{qc}) in this proposition
 is the assumption that the underlying graphs
of ${\cal A}$ and ${\cal A'}$ are equal.

By the remarks following Propositions \ref{FCR1} and \ref{FCR4}, this
proposition can be applied if ${\cal A}$ and ${\cal A}'$ are
$K$-approximations
of a  $Q$-regular $Q$-Loewner space $X$ with  $Q>1$ and 
of a $Q'$-regular space $Y$  with $Q'\le Q$, respectively.
This special case  corresponds to the
situation  in Proposition \ref{FCR0}. 
 
\proof By property (2) and (3) of a $K$-approximation, the restrictions
$p'|W$ and $p|W$ are injective. Hence $f$ is well-defined and a bijection.
 
  By Lemma \ref{weaklyup} the set $A$ is weakly $\la$-uniformly
perfect with $\la$  depending only on $s$ and  $K$.
Since $Y$ is doubling,
the subset $B$ is also doubling, quantitatively. Hence by Lemma \ref{weakqm},
in order to establish  that $f$ is uniformly \qm\, it is
enough to
show that if  $(x_1,x_2,x_3,x_4)$ is a four-tuple of distinct points
in $A$, and $[f(x_1), f(x_2), f(x_3), f(x_4)]$ is small,
 then $[x_1,x_2,x_3,x_4]$ is small,  quantitatively.
To see this let $v_i=p^{-1}(x_i)={p'}^{-1}(f(x_i))$.
 Then Proposition \ref{FCR1}  shows that if
$[f(x_1), f(x_2), f(x_3), f(x_4)]$ is  small, then
$[v_1,v_2,v_3,v_4]_Q$ is also small quantitatively. This in turn implies
by Proposition \ref{FCR1} that $[x_1,x_2,x_3,x_4]$
is small,  quantitatively.  \qed

As already mentioned,  condition (\ref{FCR11}) is true if $Q>1$ and 
 $Z$ is  $Q'$-regular  with $Q'\le Q$. This is proved in the following
proposition.   

\begin{proposition}\label{FCR2}
Suppose  $Q>1$ and let  $(Z,d,\mu)$
be  a metric measure space which is $LLC$ and
$Q'$-regular for some $Q'\le Q$.
 Let ${\cal A}$
be  a   $K$-approxi\-mation of $Z$. Then there exists
a function $\Psi\:\R^+\ra (0,\infty]$ with
$\lim_{t\to \infty} \Psi(t)=0$   depending only on $K$, $Q$ and the
data of $Z$ such that
\begin{equation} \label{FCR21}
 \GMod_Q(V_E,V_F) \le \Psi(\Delta(E,F)),
\end{equation}
whenever $E,F\sub Z$ are continua not  contained in any $K$-star.
\end{proposition}
\proof
We may assume
$\Delta(E,F)\ge 2$ and
 $R:= \diam(E)\le \diam(F)$. Fix $z_0\in E$.
Since $(G,p,r, {\cal U})$ is a $K$-approximation, we have that
\begin{equation} \label{FR2}
|d(z_0, p(u)) -d(z_0, p(v))|\le C_1  r(u) \for u,v\in V, \, u\sim v,
\end{equation}
where $C_1=C_1(K)$.
If $d(z_0, p(v)) <r(v)$ for some $v\in V$, then $U_v \cap E\ne \emptyset$,
and so $v\in V_E$. Hence $r(v)\le C_2 \diam(E)$, where $C_2=C_2(K)>0$,
because $E$ is not contained in $\St_K(v)$.
Therefore,  there exists $C_3=C_3(K)>0$ such that
\begin{equation} \label{FR3}
r(v)\le C_3 (R + d(z_0, p(v))) \for v\in V.
\end{equation}
Together with (\ref{FR2}) this shows that there exists $C_4=C_4(K) \ge 1$ such
that
\begin{equation} \label{FR4}
C_4^{-1} \le \frac{R+ d(z_0,p(v))}{R + d(z_0,p(u))}
\le C_4 \for u,v\in V, \, u\sim v.
\end{equation}
 
Now define $w:V\ra \R^+$ as follows.
Let 
$$ w(v) = \frac{ r(v) } { \log(\Delta(E,F))
(R+ d(z_0, p(v)))}
$$ if  $0\le  d(p(v),z_0) \le R\Delta(E,F)$ and let
$w(v)=0$ otherwise.
There exists $N\in \N$ such that
\begin{equation} \label{FR5}
2^{N-1} \le \Delta(E,F)< 2^N.
\end{equation}
Let $B_i:=B(z_0, 2^iR)$ for $i\in \{0,\dots,N\}$ and let $B_{-1}=\emptyset$.
 By property (2) of a $K$-approximation and by
(\ref{FR3}) there exist $C_5>0$
depending only on the data such that $U_v\sub B(z_0, C_52^iR)$ whenever
$v\in V$ and $p(v)\in B_i$.
Using (\ref{FR3}) and the $Q'$-regularity
 of $\mu$
 we obtain    for the total mass of $w$
\begin{eqnarray*}
\sum_{v\in V} w(v)^Q & \le &
  \sum_{i=0}^N \sum_{p(v)\in B_i\setminus B_{i-1}}
 w(v)^Q \\
&\lesssim & \frac{1}{(\log\Delta(E,F))^Q}
 \sum_{i=0}^N\sum_{p(v)\in B_i\setminus B_{i-1}}
\frac{r(v)^{Q'}}{(R+d(z_0,p(v)))^{Q'}}\\
&\lesssim  &\frac{1}{(\log\Delta(E,F))^Q}
 \sum_{i=0}^N\sum_{p(v)\in B_i}
\frac{\mu (U_v)}{2^{iQ'}R^{Q'}}  \\
&\lesssim  & \frac{1}{(\log\Delta(E,F))^Q}
\sum_{i=0}^N \frac{\mu(B(z_0, C_5 2^iR))}{2^{iQ'}R^{Q'}} \\
&\lesssim  & \frac{N+1}{(\log\Delta(E,F))^Q} \, \lesssim
\frac{1} {(\log\Delta(E,F))^{Q-1}}.
\end{eqnarray*}
In the last inequality we used (\ref{FR5}) and the fact $\Delta(E,F)\ge 2$.
 
On the other hand, let $u_1\sim \dots \sim u_n$ be  an arbitrary chain with
$u_1\in V_E$ and $u_n\in V_F$.
Let $d_i:=R+ d(z_0,p(u_i))$, $i\in
\{1,\dots, n\}$.
Then there is a largest number $k\in \N$ such that $d(z_0,p(u_{i}))\le
R\Delta(E,F)=\dist(E,F)$ for $i\in \{1,\dots, k\}$.
Then $d_k\gtrsim  R\Delta(E,F)$.
For otherwise, $d(z_0,p(u_k))<d_k<<R \Delta(E,F)$. This implies
 $r(u_k)\simeq \diam(U_{u_k})\gtrsim R\Delta(E,F)$ if $k=n$,
because $U_{u_k}$ then meets $F$ and contains $p(u_k)$ which is close
to $E$.
But $r(u_k)\gtrsim R\Delta(E,F)$ is also
true if $k<n$, because then by (\ref{FR2}) we have
 $r(u_k)\gtrsim |d_{k+1}-d_k| \simeq
d(z_0,p(u_{k+1})) > R \Delta(E,F)$.
 
Now $d(z_0,p(u_k))<<R \Delta(E,F)$ and $r(u_k)\gtrsim R \Delta(E,F)$ are
incompatible if $\Delta(E,F)$ is larger than a constant depending on the data,
because in this case
$E\sub N_{r(u_k)/K}(U_{u_k})\sub \St_K(u_k)$
which  is a contradiction.
 
Note that since $r(v)\lesssim \diam(E)$ for $v\in V_E$,
 we have $d_1\lesssim R$. Hence $\log(d_k/d_1)\gtrsim
\log\Delta(E,F)$, and by using (\ref{FR2}) and (\ref{FR4}) we arrive at
\begin{eqnarray*}
\sum_{i=1}^n w(v_i) & \ge & \frac{1}{\log\Delta(E,F)}\sum_{i=1}^k
\frac{r(u_i)}{d_{i}} \\
&\gtrsim & \frac{1}{\log\Delta(E,F)} \sum_{i=1}^{k-1}  \frac{|d_{i+1}-d_i|}
{d_i \wedge d_{i+1}}  \\
&\ge & \frac{1}{\log\Delta(E,F)} \sum_{i=1}^{k-1} \int_{d_i}^{d_{i+1}}
\frac{ds}s \\
&=& \frac{\log(d_k/d_1)}{\log\Delta(E,F)} \, \gtrsim\, 1.
\end{eqnarray*}
This and the mass bound for $w$ show
$$ \GMod_Q(V_E,V_F) \lesssim \frac{1}{(\log\Delta(E,F))^{Q-1}}. $$
The assertion follows from this and  $Q>1$.
\qed

In  the previous proof we used  (\ref{FR3}) in the second of the 
inequalities used to derive the  mass bound for $w$. If we do not use
(\ref{FR3}), then the proof actually shows
\begin{equation}\label{repltyson}
 \GMod_Q(V_E,V_F) \le
\left(\frac{\mesh({\cal A})} {\diam(E)\wedge \diam(F)}\right)^{Q-Q'}
\frac{C}{(\log\Delta(E,F))^{Q-1}}, 
\end{equation}
where $C$ is a constant depending only on $K$, $Q$ and the data of $Z$. 
This inequality will be useful in the proof of Theorem \ref{thmloew}. 

The goal in the proofs of Theorems \ref{mainthm} and  \ref{thmloew} is 
the construction of a quasisymmetric map between two spaces. 
Based on Proposition  \ref{unifqm} one can prove a general  result in this 
direction if one considers $K$-approximations of the spaces with 
mesh size tending to zero.

\begin{proposition} \label{qs1}
Let $Q,K,K'\ge 1$, and let   $(X,d_X)$ and $(Y,d_Y)$  be compact
 metric  measure  spaces. Assume that
${\cal A}_k=(G_k, p_k, r_k, {\cal U}_k)$ and
${\cal A'}_k=(G_k, p'_k, r'_k, {\cal U}'_k)$ for $k\in \N$ are
$K$-approximations and $K'$-approximations
 of $X$ and $Y$, respectively, whose underlying
graphs $G_k=(V^k,\sim)$ are  the same.
 
Suppose that  $X$ is connected,  and that there exist $M>0$ and some function
$\Phi$ such that $X$ and
 ${\cal A}_k$ for $k\in \N$ satisfy condition  (\ref{FCR41}). 
Suppose $Y$ is $LLC$ and doubling, and that there exist  $L>0$ and some
function $\Psi$ such that $Y$ and   ${\cal A}'_k$ for $k\in \N$ satisfy
condition (\ref{FCR11}).
 
Finally, suppose that there exists  $\la>0$ 
and vertices $v^k_1, v^k_2, v^k_3\in V^k$ for
$k\in \N$ such that  
$$ d_X(p_k(v^k_i)), p_k( v^k_j))\ge \la \diam(X) \text{ and }
d_Y(p'_k(v^k_i), p'_k(v^k_j))\ge \la \diam(Y)$$
for $k\in \N$, $i,j\in \{1,2,3\}$, $i\ne j$.
 
If $\lim_{k\to \infty}\mesh({\cal A}_k)=0$, then 
there exists an  $\eta_1$-quasisymmetric map $f\:X\ra Y$, where
$\eta_1$ depends only on the data.

If $\lim_{k\to \infty}\mesh({\cal A}'_k)=0$, then 
there exists an  $\eta_2$-quasisymmetric map $g\:Y\ra X$, where
$\eta_2$ depends only on the data.

\end{proposition}

The data here consist of  $K$, $Q$, $\la$,  the functions $\Phi$
 and $\Psi$, and the  $LLC$ and the  doubling  constant  of $Y$.
Note that we do not claim that $f$ or $g$ are surjectice.
If both $\mesh({\cal A}'_k)\to 0$ and $\mesh({\cal A}'_k)\to 0$, then 
the maps $f$ and $g$ can be constructed so that they are inverse to each other.
In this case the spaces $X$ and $Y$ are quasisymmetrically equivalent. 

The natural question arises what the relation of the conditions 
$\mesh({\cal A}_k) \to 0$  and
$\mesh({\cal A}'_k)\to 0$ is. We will later see (cf.\ 
Proposition \ref{qs2}) that even under slightly weaker assumptions
 $\mesh({\cal A}'_k)\to 0$
 actually implies $\mesh({\cal A}_k)\to 0$.   
The other direction is less clear. 

We will  apply this proposition in the case  that $X$ and $Y$ are topological
$2$-spheres. In this case $f$ and $g$ are forced to be surjective, 
since a sphere can not be embedded into a proper subset of an another
sphere of the same dimension.

\proof 
  If $\mesh({\cal A}'_k) \to 0$ or $\mesh({\cal A}'_k) \to 0$,  
 then the mutual combinatorial distance of the vertices $v^k_1, v^k_2,
v^k_3$ becomes arbitrarily large as  $k\to \infty$. 
So if $k$ is sufficiently large, $k\ge k_0$ say, then there exists  a maximal
$(2K+2L)$-separated set $W_k \sub V^k$ containing 
$v^k_1, v^k_2, v^k_3$.
Assume $k\ge k_0$ for the rest
of the proof. 

Let $ A_k:=p_k(W_k)$,  $B_k:=p'_k(W_k)$ and
$f_k\: A_k\ra B_k$, $x\mapsto p_k({p'}_k^{-1}(x))$.
Then by Proposition \ref{unifqm}, the maps  $f_k$ are $\tilde\eta_1$-\qm\
with $\tilde\eta_1$ depending on the data (and not on $k$).
Hence  the inverse maps 
$g_k=f_k^{-1}\: B_k \ra A_k$ are $\tilde\eta_2$-\qm\ with $\tilde
\eta_2$ depending
on the data.  Moreover, let  $x_i^k:=p(v_i^k)$ and $y_i^k:=p'_k(v_i^k)$ 
for $i\in \{1,2,3\}$. Then  $d_X(x_i^k, x_j^k)\ge 
\la \diam(X)$ and  $d_Y( y_i^k, y_j^k) \ge \la \diam(Y)$ for 
$i,j\in \{1,2,3\}$, $i\ne j$,  and we have  $f_k(x_i^k)=y_i^k$ and
$g_k(y_i^k)=x_i^k$. 

Every vertex $v\in V^k$ has combinatorial distance at most 
$2K+2L$ to the set $W_k$. Moreover, the sets
$U_v$, $v\in V^k, $ form a cover of $X$. It follows from the properties
of a $K$-approximation that  every point in $X$ lies within 
distance $C(K,L)\mesh({\cal A}_k)$ of the set $A_k$.
So if  $\mesh({\cal A}_k) \to 0$, then $\sup_{x\in X}\dist(x,  A_k) \to
0$ as $k\to \infty$. In this case the maps $f_k$ subconverge
to an $\tilde \eta_1$-\qm\ map  $f\:X\ra Y$ by Lemma \ref{TBA}.

Passing to appropriate 
subsequences we may assume that $x^k_i\to x_i\in X$ and $y^k_i\to y_i\in Y$
as $k\to \infty$, and $f(x_i)=y_i$ for $i\in \{1,2,3\}$. Then
$d_X(x_i, x_j) \ge \la \diam(X)$ and $d_X(y_i, y_j) \ge \la \diam(Y)$ for
$i,j\in \{1,2,3\}$, $i\ne j$.  It follows  from  remark (4) in Section
\ref{quasi-mobiusmaps} that $f$ is a $\eta_1$-quasisymmetric 
with $\eta_1$ depending on $\la$ and $\tilde\eta_1$, and hence
only on the data.

If $\mesh({\cal A}'_k) \to 0$, then by considering the maps $g_k$ one can
construct an  $\eta_2$-quasisymmetric map $g\: Y\to X$ with 
$\eta_2$ depending on the data in a similar way.  

If both $\mesh({\cal A}_k) \to 0$ and $\mesh({\cal A}'_k) \to 0$, 
then we first find a subsequence $(f_{k_l})_{l\in \N}$ of the  sequence
$f_k$ 
converging to a map $f$. Then a subsequence of the 
sequence $(g_{k_l})_{l\in \N}$
will converge to a map $g$. Then $f$ and $g$ will be 
quasisymmetries as desired, and we have in addition that
$f$ and $g$ are inverse
to each other.  \qed

\subsection{The proofs of Theorems  \ref{mainthm} and \ref{thmloew}}
\label{theproofsofthetheorems}

We will derive our Theorems \ref{mainthm} and \ref{thmloew} from more 
general theorems that give necessary and sufficient conditions for a 
metric $2$-sphere to be quasisymmetric to $\sph$.
In Theorems \ref{necsuff1} and \ref{necsuff} we will assume that 
$Z$ is linear locally connected and doubling. These conditions
are  necessary for $Z$
to be quasisymmetric to $\sph$. Moreover, a sequence of  
$K$-approximations as specified always exists under these necessary
 a priori assumptions.

\begin{theorem}
\label{necsuff1}
Let $Z$ be  metric space homeomorphic
to $\sph$ which is linearly locally connected and doubling. Suppose
$K\ge 1$ and
 ${\cal A}_k=(G_k,p_k,r_k,{\cal U}_k)$ for $k\in \N$
are  $K$-approximations
whose graphs $G_k=(V^k, \sim)$
 are combinatorially equivalent
to  $1$-skeletons of triangulations $T_k$ of $\sph$ and for which
\begin{equation} \label{mesh}
 \lim_{k\to \infty} \mesh({\cal A}_k) = 0.
\end{equation}
Suppose  there exist
 numbers $Q\ge 2$,  $k_0\in \N$, $M>0$,   and
a positive decreasing  function
$\Phi\:\R^+ \to \R^+$ 
 satisfying the following property:
 
If $k\ge k_0$  and  $E,F \sub Z$ are continua  with $\dist(V^k_E, V^k_F)\ge M$,
then
\begin{equation} \label{lowbdmod}
\Phi(\Delta(E,F)) \le  \hbox{\rm mod}^{G_k}_Q( V^k_E, V^k_F).
\end{equation}
Then there exists an $\eta$-quasisymmetric homeomorphism $f\:
Z\to \sph$ with $\eta$ depending only on the data. 

Conversely, if $Z$
is quasisymmetric to $\sph$, then condition
 (\ref{lowbdmod}) for the given sequence
${\cal A}_k$   is satisfied 
for $Q=2$, some numbers
$k_0\in \N$, $M>0$, and 
 an appropriate function $\Phi$. 
\end{theorem}

 The data in the first part of the theorem are $Q$, $K$, $M$,  $\Phi$,
and the $LLC$ and doubling constants of $Z$.

\proof 
Fix a triple $(z_1,z_2,z_3)$ of
distinct points in $Z$ such that
$d(z_i,z_j) \ge  \diam(Z)/2$ for $i,j\in \{1,2,3\}$, $i\ne j$.
Since $\mesh({\cal A}_k)\to 0$, for sufficiently
large $k$, say  $k\ge k_0$,  we can find
$v_i^k\in V^k$ such that for $x_i^k:=p_k(v_i^k)$
we have  $d(z_i, x_i^k) < \diam(Z)/4$ for
$i\in \{1,2,3\}$.
Then $d(x_i^k, x_j^k) \ge \diam(Z)/4$ for 
$i,j\in \{1,2,3\}$,  $i\ne j$.
Assume $k\ge k_0$ for the rest of the proof. 
 
The triangulation $T_k$ can be realized as a circle packing on $\sph$
(Section \ref{circlepackingsection}).
We normalize the circle packing so that  the vertices
$v_1^k, v_2^k, v_3^k$ correspond to points $y_1,y_2,y_3$ in $\sph$
 equally spaced on
some great circle.  
 The circle packings
induce canonical
$K'$-approximations ${\cal A}'_k=(G_k, p'_k, r'_k, {\cal U}'_k)$
of $\sph$, where $K'$ depends only on $K$.
Then $p'_k(v_i^k)=y_i$ and so the vertices $v^k_i$ satisfy the  condition
in Proposition \ref{qs1},  where $\la$ is a numerical constant.

Since $\sph$ is $LLC$ and $2$-regular,  and $Q\ge 2$, we see 
by Proposition \ref{FCR2} that  condition (\ref{FCR21}) is true for the space 
$\sph$ and the $K'$-approximations ${\cal A}'_k$ with $L=K'$ and a uniform 
function $\Psi$ independent of $k$. 
Therefore, the hypotheses of Proposition \ref{qs1} are satisfied for $X=Z$, 
$Y=\sph$ and our sequence of approximations. We conclude 
that there exists an  $\eta$-quasisymmetry $f\:Z \ra \sph$ where 
$\eta$ depends only on the data. Since $Z$ is a topological sphere,
this embedding has to be surjective and is hence a homeomorphism. 

Conversely, assume that there exists an
$\eta$-quasisymmetry $f\: Z \ra \sph$. Since (\ref{mesh})
implies the condition  (\ref{meshsmall}) in 
Lemma \ref{qsofapp}  for sufficiently large $k$, say for $k\ge  k_0$,
we can use the quasisymmetric images of the $K$-approximations
${\cal A}_k$ as in Lemma \ref{qsofapp}  to obtain $K'$-approximations
${\cal A}'_k=(G_k, p'_k, r'_k, {\cal U}'_k)$  of $\sph$. Here
$K'$ depends only on $K$ and  $\eta$.

Since $\sph$ is a $2$-regular $2$-Loewner space, by Proposition \ref{xferprop}
condition (\ref{FCR41}) is true for the space $\sph$ and  the
$K'$-approximations ${\cal A}'_k$
 with $Q=2$,  the  constant $M=4K'$ 
and a function $\Phi'$ independent of $k$.

Now let $k\ge k_0$,  and
 suppose that $E,F\sub Z$ are continua such that $\dist(V^k_E, V^k_E)\ge M$.
The underlying graphs of ${\cal A}_k$ and
${\cal A}'_k$ are the same. Moreover, the combinatorics of the covers
${\cal U}_k$ and ${\cal U}'_k$ correspond under the mapping $f$.
This shows that for $E'=f(E)$ and $F'=f(F)$ we have
$V^k_E=V^k_{E'}$, $V^k_F=V^k_{F'}$,  and $\dist(V^k_E, V^k_F)=
\dist(V^k_{E'}, V^k_{F'})\ge M,$ where the sets $V^k_E$ etc.\ are interpreted
with respect to the appropriate approximations. 
Hence we get 
$$\Phi'(\Delta(E',F'))\le \hbox{\rm mod}^{G_k}_2(V^k_{E'}, V^k_{F'})=
\hbox{\rm mod}^{G_k}_2(V^k_{E}, V^k_{F}).
$$
Condition  (\ref{lowbdmod})  for an appropriate function $\Phi$ independent 
of $k$ will follow from this,
if we can show that $\Delta(E,F)$ is large if and only if $\Delta(E', F')$
is large, quantitatively.  
But this last statement follows from the quasisymmetry of $f$ and Lemma
\ref{sepacont}. 
 \qed

\medskip
 As an immediate application of this theorem we get a proof of Theorem
\ref{thmloew}.

{\em Proof of Theorem \ref{thmloew}.\ } 
Suppose  $Z$  is $Q$-regular and $Q$-Loewner for $Q\ge 2$. 
Then  $Z$ is $LLC$ and  doubling.
Corollary \ref{goodapprox} shows that  there exists  $K\ge 1$ and a  sequence
 of $K$-approximations ${\cal A}_k=(G_k, p_k, r_k, {\cal U}_k)$
 whose  graphs $G_k=(V_k,\sim)$ are combinatorially equivalent to 
$1$-skeletons of triangulations $T_k$ of $Z$ and for which
(\ref{mesh}) is true. Now the $Q$-regularity of  $Z$, Proposition 
\ref{xferprop} and the $Q$-regularity  show that condition 
 (\ref{FCR41}) is true
for the  $K$-approximations ${\cal A}_k$ with $M=4K$ and a 
function $\Phi$ independent of $k$.
 Theorem \ref{necsuff1} implies that there exists a quasisymmetric
homeomorphism $f\: Z\ra \sph$.
 A  result
by Tyson \cite{tys} shows that if
a $Q$-regular $Q$-Loewner space is quasisymmetrically mapped onto a
$Q'$-regular space, then $Q'\ge Q$. But $\sph$ is $2$-regular,
and so we can apply this for $Q'=2$ and get $2\ge Q$.  
 Since also $Q\ge Q'=2$  by assumption, 
 we must have $Q=2$.
The proof of Theorem \ref{thmloew} is complete.
\qed 

It may be worthwhile to point out that in the previous proof an argument
can be  given  that avoids 
invoking Tyson's theorem.

 Suppose $Z$ is $Q$-regular $Q$-Loewner space
and $f\:Z\ra \sph$ a quasisymmetric homeomorphism. 
Let ${\cal A}_k$ be a sequence of $K$-approximations  of $Z$ with 
underlying graphs $G_k=(V^k,\sim)$  such that 
$\lim_{k\to \infty} \mesh({\cal A}_k)=0$. Let ${\cal A}'_k$ be 
the $K'$-approximation of $\sph$ obtained as the image of ${\cal A}_k$ 
under $f$. Then $\lim_{k\to \infty} \mesh({\cal A}'_k)=0$.
Let $E,F\sub Z$ be two disjoint continua and $E':=f(E)$, $F':=f(E)$.
Then  by Proposition  \ref{xferprop} and by the remark following the proof 
of Proposition \ref{FCR2} we have for sufficiently large $k$
\begin{eqnarray*}
  \Phi(\Delta(E,F))& \le & \Mod_Q(E,F) \lesssim
\GMod^{G_k}_Q(V^k_E, V^k_F)=  
\GMod^{G_k}_Q(V^k_{E'}, V^k_{F'}) \\ & \lesssim &  
\left(\frac{\mesh({\cal A}'_k)} {\diam(E')\wedge \diam(F')}\right)^{Q-2}
\frac{1}{(\log\Delta(E',F'))^{Q-1}}.
\end{eqnarray*}
Here $\Phi$ is  a positive function  provided by the 
$Q$-Loewner property of $Z$.
Moreover, the multiplicative constants implicit in this
inequality are independent of $E$, $F$ and $k$.
Note that the additional assumptions 
in Propositions \ref{xferprop} and \ref{FCR2}
are true for our  continua
if $k$ is sufficiently large. If $Q>2$ then the  last term 
in the inequality tends to zero, since the mesh size tends to zero.
But this is impossible, since the first term is independent of $k$ and 
positive. Hence $Q=2$.

\begin{theorem}
\label{necsuff}
Let $Z$ be  metric space homeomorphic
to $\sph$ which is linearly locally connected and doubling. Suppose
$K'\ge 1$,  and
 ${\cal A}_k=(G_k,p_k,r_k,{\cal U}_k)$ for $k\in \N$
are   $K$-approximations
whose graphs  $G_k=(V^k,\sim)$ are  combinatorially equivalent to 
 the  $1$-skeletons  of  triangulations $T_k$ of $\sph$ and for which
\begin{equation} \label{mesh1}
 \lim_{k\to \infty} \mesh({\cal A}_k) = 0.
\end{equation}
 
Suppose that there exist 
numbers $k_0\in \N$, $L>0$, and
a  function $\Psi\:\R^+ \to (0,\infty]$ with $\lim_{t\to \infty}
\Psi(t)=0$ satisfying the following property:
 
If $k\ge k_0$  and  $E,F \sub Z$ are continua not contained
in any $L$-star of ${\cal A}_k$, then
\begin{equation} \label{uppbdmod}
 \hbox{\rm mod}^{G_k}_2( V^k_E, V^k_F) \le \Psi(\Delta(E,F)).
\end{equation}

Then there exists an $\eta$-quasisymmetric homeomorphism $g\: Z\to \sph$
with $\eta$ depending only on the data. 

Conversely, if $Z$
is quasisymmetric to $\sph$, then condition
 (\ref{uppbdmod}) for the given sequence
${\cal A}_k$   is satisfied
for some numbers
$k_0\in \N$, $L>0$, and
 an appropriate function $\Psi$.
\end{theorem}

 The data in the first part of the theorem are  $K$, $L$,  $\Psi$,
and the $LLC$ and doubling constants of $Z$.

\proof  The proof of 
 this theorem is very similar to the proof 
of Theorem \ref{necsuff1}.  
For the sufficency part note again 
that  the triangulation $T_k$ can be realized as a normalized 
 circle packing on $\sph$.
The circle packings
induce canonical
$K'$-approximations ${\cal A}'_k=(G_k, p'_k, r'_k, {\cal U}'_k)$
of $\sph$, where $K'$ depends only on $K$.

As in the proof of Theorem \ref{necsuff1},  for sufficiently large $k$
 we can  find 
vertices
$v_1^k, v_2^k, v_3^k \in V^k$ satisfying the 
 the condition  in Proposition \ref{qs1} where  $\la>0$ is a numerical
constant.  
Since $\sph$ is  $2$-regular and $2$-Loewner, Proposition \ref{xferprop}
implies 
that condition  (\ref{FCR41}) is true for the space
$\sph$ and the $K'$-approximations ${\cal A}'_k$ with  $M=4K'$ 
and a function $\Phi$ independent of $k$.

It follows that the hypotheses of Proposition \ref{qs1} are satisfied
for  $X=\sph$ and the $K'$-approximations  ${\cal A}'_k$
and $Y=Z$ and the  $K$-approximations ${\cal A}_k$. (Note that the roles
of ${\cal A}$ and  ${\cal A}'_k$ in this proof and in Proposition \ref{qs1}
are reversed). Since $\mesh({\cal A}_k)\to 0$ it follows 
that there exists an  $\eta$-quasisymmetry $g\:Z \ra \sph$ where
$\eta$ depends only on the data.  Again $g$ has to be 
a homeomorphism.

For the  converse assume that
 there exists an
$\eta$-quasisymmetry $g\: Z \ra \sph$. Again 
for sufficiently large $k$
we  obtain  $K'$-approximations ${\cal A}'_k$ of $\sph$
 with
$K'=K'(\eta,K)$ as the  quasisymmetric images under $g$ of the 
$K$-approximations ${\cal A}_k$. 
The sphere $\sph$ is $2$-regular, so by Proposition 
\ref{FCR4} we have
condition (\ref{FCR41}) for $Q=2$, $L:=K'$ and an appropriate 
function $\Psi'$ independent of $k$.
Now suppose $E,F$ are continua not contained in any $L$-star
with respect to ${\cal A}_k$. We have 
 $\St_{L}(v)=g(\St_{L}(v))$, where the first star is with respect
to ${\cal A}'_k$ and hence a subset of $\sph$, and the second star 
is with respect
to ${\cal A}_k$ and hence a subset of $Z$.
This implies that $E'=g(E)$ and $F'=g(F)$ are not contained in any 
$L$-star with respect to ${\cal A}'_k$.   
Hence  
$$
\hbox{\rm mod}^{G_k}_2( V^k_E, V^k_F)=
\hbox{\rm mod}^{G_k}_2( V^k_{E'}, V^k_{F'}) \le \Psi(\Delta(E',F')).
$$
Now $\Delta(E',F')$ is large if and only if $\Delta(E,F)$, quantitatively.
Hence  condition (\ref{uppbdmod})
follows with $L=K'$, and an appropriate 
function $\Phi$ independent of $k$. \qed 

\medskip
{\em Proof of Theorem \ref{mainthm}.\ }
As we remarked in the introduction,  only the sufficiency part of  Theorem 
\ref{mainthm}
demands a proof. Since linear local contractibility and 
linear local connectivity are quantitatively equivalent for 
topological $2$-spheres, we can assume that $Z$ is $LLC$.
We will show that there exists an $\eta$-quasisymmetric
homeomorphism $g\: \Z\to \sph$, where $\eta$ depends only on
the data. Here we call
  the $LLC$ constant,
 and the constant that enters the condition for $2$-regularity
(where $\mu={\cal H}^2)$ the data of $Z$.

Note that $Z$ is doubling with a constant only 
depending on the data.  
Corollary \ref{goodapprox} shows that  there exists  $K\ge 1$ depending
on the data 
 and a  sequence  
 of $K$-approximations ${\cal A}_k=(G_k, p_k, r_k, {\cal U}_k)$
 whose  graphs $G_k=(V^k,\sim)$ are
$1$-skeletons of triangulations $T_k$ of $Z$ and for which
(\ref{mesh1}) is true.   Since $Z$ is $LLC$ and $2$-regular, 
Proposition \ref{FCR2} shows that the condition (\ref{uppbdmod}) 
is true for $L=K$ and an appropriate function $\Phi$ depending on the data.
Now Theorem \ref{necsuff} shows that there exists a $\eta$-quasisymmetric 
homeomorphism $g\: \Z\to \sph$, where $\eta$ depends only on 
the data.   \qed

Theorem \ref{mainthm} is quantitative
as the  proof above  shows.
Namely, if $Z$ is a metric 
space homeomorphic to $\sph$ that  is Ahlfors $2$-regular and $LLC$, then 
there exists an $\eta$-quasisymmetric homeomorphism $g\: Z\to \sph$, where
$\eta$ depends only on the data, i.e., the 
 constants in the Ahlfors $2$-regularity and the $LLC$ conditions.
Conversely,  if $Z$ is a metric space 
for which  there exists an $\eta$-quasisymmetric homeomorphism
$g\: \Z\to \sph$, then $Z$ is $\la$-$LLC$ with $\la$ only depending on $\eta$.

\subsection{Asymptotic conditions}  
\label{asymptoticconditions}

Cannon's  paper \cite{cannonacta} provides a
framework that allows one to speak of  modulus for
subsets of a topological space.  A shingling
${\cal S}$ of a topological space $X$ is a locally
finite cover consisting of compact connected subsets of $X$.
When $X=\sph$ and  $R\subset \sph$
is an annulus,  Cannon defines invariants
$M(S,R)$ and $m(S,R)$ which are combinatorial analogs
for the classical moduli of annuli.
He then studies a sequence of shinglings $S_j$ of  $\sph$
with mesh size tending to zero.  His main theorem---the
combinatorial Riemann mapping theorem---is
a necessary and sufficient condition for the existence
of a homeomorphism $f:\sph\ra\sph$ such that for every
annulus $R\subset \sph$, the moduli
$M(f_*(S_j,R))$ and $m(f_*(S_j,R))$
agree with the standard $2$-modulus
to within a fixed multiplicative factor, for sufficiently
large $j$.

The combinatorial Riemann mapping theorem is similar in
spirit to Theorems \ref{necsuff1} and \ref{necsuff}: all three
results give necessary and sufficient conditions for
a ``conformally flavored'' structure on the $2$-sphere to be equivalent
modulo a homeomorphism to the standard structure.
Any of these  theorems can be used to give
necessary and sufficient conditions for a Gromov hyperbolic
group to admit a discrete,
cocompact, and isometric action on hyperbolic space
$\H^3$.   The paper \cite{canswe} uses  \cite{cannonacta}
and \cite[Corollary, p. 468]{sul} to  give such conditions; the conditions
in \cite{canswe} are in turn applied in \cite{suffrich}.
Our Theorems \ref{necsuff1} or \ref{necsuff} can be combined
directly  with Sullivan's theorem.   The point here is that the action
$G\acts\geo G$ of a non-elementary hyperbolic group on its
boundary is by uniformly \qm\ homeomorphisms, and if one conjugates
this action by a \qs\ homeomorphism $\geo G\ra\sph$, the resulting
action $G\acts\sph$ is also uniformly \qm, in particular
uniformly quasiconformal, so that \cite{sul} may be applied.

On the other hand, there are significant differences between our
approach and
Cannon's approach. Cannon's 
hypotheses and conclusions do not involve metric information, and only
relate to the  limiting behavior of the combinatorial
moduli.  In contrast, Theorems \ref{necsuff1} and \ref{necsuff}
hypothesize inequalities between relative distance
(which is metric based) and combinatorial modulus
which hold  uniformly for every $K$-approximation
in the given sequence; and they assert that
the metric space is \qs\ to $\sph$, which is a metric conclusion.
 
The interesting parts of Theorems \ref{necsuff1} and \ref{necsuff}  are 
the sufficient conditions. An upper bound  for a modulus is easier to establish
than a lower bound, because for a lower bound an inequality for the total
mass of  {\em all} admissible test functions has to be shown whereas
an upper bound  already follows from a mass bound for {\em one} 
test function. In this respect, Theorem \ref{necsuff} seems to be more 
useful, because its hypotheses require  upper  
modulus bounds. In view of Cannon's work it seems worthwhile to find
a sufficient condition in the spirit of Theorem \ref{necsuff} that works with 
an asymptotic condition for the graph modulus as in (\ref{uppbdmod}).
The following theorem provides
such a result where we further weaken  the requirements 
for which sets $E$ and $F$ an asymptotic modulus inequality has to hold. 

\begin{theorem}
\label{suff}
Let $Z$ be a  metric space homeomorphic
to $\sph$ which is linearly locally connected and doubling. Suppose
$K\ge 1$,  and
 ${\cal A}_k=(G_k,p_k,r_k,{\cal U}_k)$ for $k\in \N$
are   $K$-approximations
whose graphs  $G_k=(V^k,\sim)$ are  combinatorially equivalent to
 the  $1$-skeletons  of  triangulations $T_k$ of $\sph$ and for which
\begin{equation} 
 \lim_{k\to \infty} \mesh({\cal A}_k) = 0.
\end{equation}

Suppose there exist a   numbers  $C>0$ and $\la >1$  
with the following property:
If $B=B(a,r)$ and $\la B=B(a,\la r)$ are balls in $Z$, 
then we have
\begin{equation} \label{uppbdmod2}
 \limsup_{k\to \infty}\, 
 {\rm mod}^{G_k}_2( V^k_B, V^k_{Z\setminus \la B}) <C. 
\end{equation}
Then there exists an $\eta$-quasisymmetric homeomorphism $g\: Z\to \sph$
with $\eta$ depending only on the data.

Conversely, 
if $Z$
is quasisymmetric to $\sph$, then there exist $C>0$ and $\la >1$ such that
 condition
 (\ref{uppbdmod2}) is satisfied for the given sequence
${\cal A}_k$.
\end{theorem}

\no
The data are $K$, $C$, $\la$,   the $LLC$ constant, and  the doubling constant.

 If  $B$ is a ball in $Z$,  let $A$ be the ``annulus" $A=\la B\setminus B$.
Its complement consists of the disjoint sets 
$B$ and $Z\setminus \la B$. The  $2$-modulus of the curve family $\Gamma$ 
joining $B$ and $Z\setminus \la B$ can be considered as 
 the $2$-modulus of the annulus $A$. The appropriate combinatorial 
version of this modulus with respect to the $K$-approximation 
${\cal A}_k$  is ${\rm mod}^{G_k}_2( V^k_B,
V^k_{Z\setminus \la B})$ which appears in (\ref{uppbdmod2}). 
So this inequality essentially says that the combinatorial analog of 
the $2$-modulus of $A$ is asymptotically  bounded above by a fixed constant. 

In order to prove this theorem we have to revisit some of the 
material in Section \ref{theferrandcross-ratio} and
prove asymptotic versions.    

\begin{proposition}\label{FCR1A}
Let   $Z$
be    a locally compact  metric measure
  space which is $\la$-$LLC$, $\la\ge 1$.
Suppose $K\ge 1$,  and 
 ${\cal A}_k=(G_k, p_k, r_k, {\cal U}_k)$ for $k\in \N$ are 
$K$-approximations  of $Z$ with graphs  $G_k=(V^k, \sim)$. Assume
that $\mesh({\cal A}_k)\to 0$ as $k\to \infty$.

 Let $Q\ge 1$, and suppose
that there exists  function $\Psi\:\R^+\ra (0,\infty]$ with
$\lim_{t\to \infty} \Psi(t)=0$ such that
\begin{equation}\label{FCR11A}
\limsup_{k\to \infty}\, {\rm mod}_Q^{G_k}(V^k_E,V^k_F) \le \Psi(\Delta(E,F)),
\end{equation}
whenever $E,F\sub Z$ are disjoint  continua.
 
 Then there exists
a function $\phi\:\R_0^+\ra [0,\infty]$ with $\lim_{t\to 0} \phi(t)=\phi(0)=0$ 
depending only on $K$, $Q$, $\Psi$ and the
data of $Z$  with the following property:
 
Suppose $(z_1, z_2, z_3, z_4)$ is a four-tuple of 
 points in $Z$ with $\{z_1,z_3\} \cap \{z_2,z_4\}= \emptyset$,
and  assume  that for $k\in \N$ 
and $i\in \{1,2,3,4\}$  we have vertices $v^k_i\in V^k$ such 
that $ p_k(v^k_i) \to  z_i$ for $k\to \infty$, $i\in \{1,2,3,4\}$.
Then
$$
\limsup_{k\to \infty}\, [v_1^k, v_2^k,v_2^k,v_3^k]^{G_k}_Q \le
\phi([z_1, z_2, z_3, z_4]). 
$$
\end{proposition}
 
We want to allow the possibility $z_1=z_3$ or $z_2=z_4$ here. In this case
we set $[z_1, z_2, z_3, z_4]=0$, which is a consistent extension of the 
definition of the cross-ratio. 
Note that $[v_1^k, v_2^k,v_2^k,v_3^k]^{G_k}_Q$ is a cross-ratio with respect 
to $G_k$. 
The proposition says that if $[z_1, z_2, z_3, z_4]$ is small, then 
$[v_1^k, v_2^k,v_2^k,v_3^k]^{G_k}_Q$ is asymptotically small, quantitatively.

\proof
If $[z_1, z_2, z_3, z_4]$ is small, then by Lemma \ref{goodcont}
there exist continua $E'$ and $F'$ with $z_1,z_3\in E$, $z_2,z_4\in F$
and $\Delta(E',F')$ large, quantitatively. If $z_1=z_3$ or $z_2=z_4$
then $\Delta(E',F')$ can be made arbitrarily large. 
Since $Z$ is locally compact and  $LLC$ and hence locally connected,
 we can find compact connected neighborhoods   $E$ and 
$F$ of $E'$ and $F'$, respectively, such that $\Delta(E,F)$ is large,
quantitatively.
Since $\mesh({\cal A}_k)\to 0$ we will have $p_k(v^k_1)\in
U_{v_1^k}\cap E$
and $p_k(v^k_3)\in U_{v_3^k}\cap E$ for large $k$. 
In particular, $v_1^k, v_3^k\in V^k_E$. Similarly, $v_2^k, v_4^k\in V^k_F$ 
for large $k$.
The rest of the proof now proceeds as the proof of Proposition \ref{FCR1}.
For large $k$ we can find chains $A_k\sub N_K(V^k_E)$ connecting 
$v_1^k, v^k_3$ and chains 
$B_k\sub N_K(V^k_F)$ connecting 
$v_2^k, v^k_4$. Then 
 by Lemma \ref{setchange} we 
have $$[v^k_1,v^k_2,v^k_3,v^k_4]^{G_k}_Q
\le {\rm mod}^{G_k}_Q(A_k, B_k)\le C(K)
{\rm mod}^{G_k}_Q(V^k_E, V^k_F).$$  So  our assumptions imply
$$\limsup_{k\to \infty}\, [v_1^k, v_2^k,v_2^k,v_3^k]^{G_k}_Q \le  C(K)
\Psi(\Delta(E,F))
.$$
Since $\Delta(E,F)$ is large and
 $\Psi(t)\to 0$ as $t\to \infty$ we get the desired quantitative 
conclusion.
\qed

The following proposition corresponds to one of the parts  of
 Proposition \ref{qs1}.
We have replaced  condition (\ref{FCR11})  by the asymptotic condition
(\ref{FCR11A}). 

\begin{proposition} \label{qs2}
Let $Q,K\ge 1$, and let   $(X,d_Y)$ and $(Y,d_Y)$  be compact
 metric  measure  spaces. Assume that
${\cal A}_k=(G_k, p_k, r_k, {\cal U}_k)$ and
${\cal A'}_k=(G_k, p'_k, r'_k, {\cal U}'_k)$ for $k\in \N$ are
$K$-approximations of $X$ and $Y$, respectively, whose underlying
graphs $G_k=(V^k,\sim)$ are  the same. Moreover, assume
$\lim_{k\to \infty}\mesh({\cal A}'_k)=0.$
 
Suppose $X$ is connected,  and there exists $M>0$ and some function $\Phi$
such that  $X$ and
 ${\cal A}_k$ for $k\in \N$  satisfy condition  (\ref{FCR41}).
 Suppose $Y$ is $LLC$ and doubling,
  and  $Y$ and ${\cal A}'_k$ satisfy
condition (\ref{FCR11A})  for  some  function $\Psi$.
 
Suppose that there are vertices $v^k_1, v^k_2, v^k_3\in V^k$ for
$k\in \N$ such that  that for some constant $\la>0$ we have
$$ d_X(p_k(v^k_i)), p_k( v^k_j))\ge \la \diam(X) \text{ and }
d_Y(p'_k(v^k_i), p'_k(v^k_j))\ge \la \diam(Y)$$
for $k\in \N$, $i,j\in \{1,2,3\}$, $i\ne j$.
 
If $\lim_{k\to \infty}\mesh({\cal A}'_k)=0$, 
then there exists a $\eta$-quasisymmetric map $f\:X\ra Y$, where
$\eta$ depends only on the data.
\end{proposition}
 
The data here consist of  $K$, $Q$, $\la$,  the functions $\Phi$
 and $\Psi$, and the  $LLC$ and the  doubling  constant  of $Y$.

In the proof we will show that $\mesh({\cal A}_k)\to 0$. Since 
condition (\ref{FCR11}) is stronger than condition (\ref{FCR11A}),
this justifies the remark after Proposition \ref{qs1}. Namely, that
that under the assumptions of this  proposition
we have that $\mesh({\cal A}'_k)\to 0$ implies $\mesh({\cal A}_k)\to 0$. 
 
\proof
1. In this proof we will call distortion functions  those functions 
   $\phi\:\R^+_0\ra  [0,\infty]$
for which $\phi(t)\to \phi(0)=0$ as $t\to 0$.
   We will first establish the existence of a distortion 
function $\phi_1$ depending on the data
 with the following property.
If $u_1^k,u_3^k\in V^k$ for $k\in \N$,  and $p_k(u^k_i)\to z_i$
and $p_k'(u^k_i)\to w_i$
as $k\to \infty$ for $i\in\{1,3\}$,
then
\begin{equation} \label{MainEst}
\frac{ d_X(z_1, z_3)}{
\diam(X)} \le \phi_1\left(\frac{ d_Y(w_1,w_3)}{\diam(Y)}  \right).
\end{equation}
 
To prove this we may assume $d_Y(w_1,w_3)< (\la/3) \diam(Y)$.
Hence if $w^k_i:=p'_k(u^k_i)$ for $i\in \{1,3\}$
 we have $d_Y(w^k_1,w^k_3)< (\la/3) \diam(Y)$
for large $k$. For such $k$ there will be at least two among
the  vertices $ v^k_1, v^k_2, v^k_3$, call them $u_2^k$ and $u_4^k$,
 such that
we have $\dist(\{w^k_1,w^k_3\}, \{w^k_2,w^k_4\}) \ge (\la/3) \diam(Y)$,
where
we set $w_i^k=p'_k(u^k_i)$ also for $i\in \{2,4\}$.
Then for large $k$ we obtain
$$ [w^k_1,w^k_2,w^k_3,w^k_4]\le C(\la) \frac{ d_Y(w^k_1,w^k_3)} {\diam(Y)}.$$
We may assume that we have limits  $w^k_2 \to w_2$ and $w^k_4\to w_4$
for $k\to \infty$. Then $\{w_1,w_3\} \cap \{w_2,  w_4\}=\emptyset$, 
 and so 
Proposition \ref{FCR1A} and the previous inequality show that there exist
a distortion functions  $\phi_2$ and  $\phi_3$ depending on the data
 such 
$$\limsup_{k\to \infty}\,
[u^k_1,u^k_2,u^k_3,u^k_4]^{G_k}_Q
\le \phi_2([w_1,w_2,w_3,w_4]) \le  \phi_3\left
(\frac{ d_Y(w_1,w_3)}{\diam(Y)}  \right).$$
  
Since $\mesh({\cal A}'_k)\to 0$ as $k\to \infty$ and the points
$w_1, w_2, w_4$ are distinct, 
the combinatorial separation of the vertices
 $u_1^k, u_2^k, u_4^k$ becomes arbitrarily
large as $k\to \infty$. We make the momentary extra assumption that
the combinatorial separation of $u_1^k$ and $u_3^k$ is at least
$4K$. Let $z_i^k=p_k(u^k_i)$. Note that $d_X(z_2^k, z_4^k) \ge (\la/2) \diam(X)$
for large $k$ by choice of $u_2^k$ and $u_4^k$.
Then from   Proposition \ref{FCR4}  we infer  that for
sufficiently large $k$
$$ \frac{ d_X( z^k_1 , z^k_3)}{\diam(X)} \le
C(\la) [z^k_1 ,z_2^k, z^k_3, z^k_4] \le
\phi_4([u^k_1,u^k_2,u^k_3,u^k_4]^{G_k}_Q),
$$
where $\phi_4$ is a distortion function depending on the data. 
Letting  $k$ tend to infinity,  the claim (\ref{MainEst})
 follows under the additional assumption on the
combinatorial separation of $u_1^k$ and $u_3^k$.
 
2. In order to establish the general case of
(\ref{MainEst}), we first show that $\mesh({\cal A}_k) \to 0$ as
$k\to \infty$.
Arguing by contradiction
and passing to a subsequence if necessary, we may assume there there exists
$\delta>0$ and $a_1^k\in V^k$ with $r_k(a_1^k)\ge \delta>0$ for $k\in \N$.
 Since the mesh size of ${\cal A}'_k$ tends to $0$,
the cardinality of $G_k$ tends to infinity. Moreover,
 $G_k$ is connected
and its vertex degree is bounded.
Thus, for sufficiently large
$k$ we can find a vertex $a_3^k\in V^k$ with $4K\le
 k_{G_k}(a^k_1, a_3^k)\le
5K$.
 Then $U_{a^k_1}\cap U_{a_2^k} =\emptyset$ and it follows
$ d_X(p_k(a_1^k), p_k(a_3^k))\ge r_k(a_1^k)\ge \delta.$
Letting $x_i^k:=p_k(a_i^k)$  and
$y_i^k:=p'_k(a_i^k)$ and passing to subsequences,
we may assume that $x_i^k\to x_i$ and  $y_i^k\to y_i$ for $k\to \infty$,
$i\in\{1,3\}$.  Then $d_X(x_1, x_3)\ge \delta>0$. On the other hand,
$y_1=y_3$, since the combinatorial distance of $a^k_1$ and
$a_3^k$ is uniformly
bounded by choice of $a_3^k$,
 and the mesh size of ${\cal A}'_k$ tends to zero.
But the combinatorial distance of $u^k_1$ and $u_3^k$ was at least $4K$ for
large $k$, so
we can apply (\ref{MainEst}) and get a contradiction.
 
3. Once we know that the mesh size of ${\cal A}_k$ tends to zero, we can
verify (\ref{MainEst})  without the additional assumption on
the combinatorial separation of $u_1^k$ and $u_3^k$.
For if $z_1=z_3$, then there is nothing
to prove. If $z_1\ne z_3$, then the combinatorial distance of $u^k_1$
and $u_3^k$
becomes arbitrarily large, since $\mesh({\cal A}_k)\to 0$ as $k\to \infty$.

4. Let $A$ be a countable dense subset of $X$.
For $z\in A$ and $k\in \N$ we can find $u_k(z)\in V^k$ with
$z\in U_{u_k(z)}$. Since $\mesh({\cal A}_k)\to 0$, we have
 $p_k(u_k(z))\to z$ as $k\to \infty$,   $z\in A$.
  Define $f_k(z):= p'_k(u_k(z))$.
By passing to successive subsequences and taking a final   ``diagonal 
subsequence'' we may assume that  the countably 
many sequences $(f_k(z))_{k\in \N}$, $z\in A$,
 converge, $f_k(z)\to f(z)$ say, as $k\to
\infty$. From (\ref{MainEst}) and  the definition of $f$, we get 
 (\ref{MainEst}) for arbitrary $z_1,z_3\in A$
and $w_1=f(z_1)$ and $w_3=f(z_3)$. In particular,  
$f\: A\to Y$ is injective.

5. We claim that the map $f$ is
$\tilde\eta$-\qm\ with $\tilde \eta$ only depending on the data.
To see this note that
the  set $A$ is weakly $\la'$-uniformly
perfect with a  fixed constant,  $\la'=2$ say.
Since $Y$ is doubling,
the subset $f(A)$ is also doubling, quantitatively. Hence by Lemma \ref{weakqm},
in order to establish  that $f$ is uniformly \qm\, it is
enough to
show that if  $(x_1,x_2,x_3,x_4)$ is a four-tuple of distinct points
in $A$, and $[f(x_1), f(x_2), f(x_3), f(x_4)]$ is small,
 then $[x_1,x_2,x_3,x_4]$ is small,  quantitatively.
By definition of $f$, we
we can find  $u_i^k\in V^k$ such $x_i\in U_{u_i^k}$
and $p'_k(u^k_i) \to y_i:=f(x_i)$ for
$k\to \infty$, $i\in \{1,\dots, 4\}$.
 Then Proposition \ref{FCR1A}  shows that if
$[y_1, y_2, y_3, y_4]$ is  small, then
$\limsup_{k\to\infty} [u^k_1,u^k_2,u^k_3,u^k_4]^{G_k}_Q$ is also small,
 quantitatively. Since the points $y_i$ are distinct, the
combinatorial separation of the vertices $u^k_i$ is arbitrarily large
for $k\to \infty$.   This implies
by Proposition \ref{FCR4} that  $[p_k(u^k_1), p_k(u^k_2), p_k(u^k_3),
p_k(u^k_4)]$ for large $k$ is small, quantitatively.
Passing to the limit we conclude that 
$$[x_1,x_2,x_3,x_4]=\lim_{k\to \infty}\, 
[p_k(u^k_1), p_k(u^k_2), p_k(u^k_3), p_k(u^k_4)]$$
is small,  quantitatively.

6. There are points $z_1,z_2,z_3$ in $A$ whose mutual distance is at
least $\diam(X)/4$.
The estimate (\ref{MainEst}) and the definition of $g$  show that
the mutual distance of the points $f_i(z_1), f_i(z_2), f_i(z_3)$
is bounded below by $c\diam(Y)$, where $c>0$ is a constant depending
on the data. Hence $f\: A\ra Y$ is $\eta$-quasisymmetric 
with $\eta$ depending on the data. 
Since $A$ is dense and $Y$ is compact, there is a unique extension 
of $f$ to an $\eta$-quasisymmetric  map on $X$.
Calling this map also $f$, we get  the desired
 quasisymmetry.
 \qed

\medskip
\no {\em Proof of Theorem \ref{suff}.\ }
To prove sufficiency, we want to apply Proposition \ref{qs2} for
$Q=2$, $X=\sph$ and $Y=Z$. As in the proof of Theorem \ref{necsuff1} one can 
realize the triangulations $T_k$ as normalized circle packings.
 The circle packings
induce canonical
$K'$-approximations ${\cal A}'_k=(G_k, p'_k, r'_k, {\cal U}'_k)$
of $\sph$, where $K'$ depends only on $K$. 
Again as in the proof of Theorem \ref{necsuff1} we can use suitable
normalizations so that for sufficiently large $k$
we can find vertices
$v_1^k, v_2^k, v_3^k \in V^k$ 
 satisfying the
 the condition  in Proposition \ref{qs2} where  $\la>0$ is a numerical
constant.
Since $\sph$ is  $2$-regular and $2$-Loewner, Proposition \ref{xferprop}
implies
that condition  (\ref{FCR41}) is true for the space
$X=\sph$ and the $K'$-approximations ${\cal A}'_k$ with  $M=4K'$
and a function $\Phi$ independent of $k$.

Since $\mesh({\cal A}_k)\to 0$ the only thing that remains to
be verified is  that with $Y=Z$, the $K$-approximations ${\cal A}_k$
satisfy the asymptotic condition (\ref{FCR11A})  for some function 
$\Psi$ depending on the data. 

To see that this is true, let  $E$ and $F$ be arbitrary disjoint 
 continua.
We have to show 
that the combinatorial modulus ${\rm mod}_2^{G_k}(V^k_E,V^k_E)$
for large $k$  is small  if the 
relative distance of $E$ and $F$ is large, quantitatively. 

We may assume $\diam(E)\le \diam(F)$. 
Pick $a\in E$,  let $r=2\diam(E)$ and  $B_i:= B(a,  \la^{2i-2}r)$
for $i\in \N$.
Then $E\sub B_1$ and  $B_i\sub \la B_i \sub  \la^2 B_i = B_{i+1}$ for 
$i\in \N$. Let $N$ be the largest integer such that $r \la^{2N-1}<
\dist(E,F)$. Note that $N$ is large if and only if  $\Delta(E,F)$ is large,
quantitatively. Then
$$E\sub B_1\sub \la B_1  \sub \la^2B_1= B_2 \sub \la B_2 
\sub  
 \dots \sub  B_N \sub 
\la B_N \sub  Z\setminus  F.$$

Since $\mesh({\cal A}_k)\to  0$, there exists $k_1\in \N$ such that 
if $k\ge k_1$ and $v\in V^k_{\la B_i}$ for some $i\in  \{1,\dots, N-1\}$,
then $v\notin V^k_{Z\setminus B_{i+1}}$.
For suppose $v\in V^k_{\la B_i} \cap V_{Z\setminus B_{i+1}}$.
Then $U_v\cap \la B_i \ne \emptyset$ and
 $U_v \cap (Z\setminus B_{i+1}) \ne \emptyset$.
Hence
$2K r_v \ge \diam(U_v)\ge \la^{2i}(1-1/\la) r \ge (1-1/\la) r.$
This is impossible if $\mesh({\cal A}_k)$ is small enough. 

By our assumption we can find $k_2\in \N$ such that  for $k\ge k_2$
and $i\in \{1,\dots, N\}$
we have
$\text{mod}^{G_k}_2(V_{B_i}, V_{Z\setminus \la B_i})< C$. 
Consider a fixed $K$-approximation ${\cal A}_k$ for  $k\ge k_3:= k_1\vee k_2$.
To simplify notation we drop the sub- or superscript $k$.

By choice of $k$, there exists  a weight $w_i\: V\to [0,\infty)$ which is 
admissible for the pair $(V_{B_i}, V_{Z\setminus \la B_i})$ and satisfies
$$ \sum_{v\in V} w_i(v)^2<C.$$
Define $w(v):=\sup_{i \in \{1,\dots, N\}} w_i(v)$ for $v\in V$. 
Then 
\begin{equation}\label{mbd}
\sum_{v\in V} w(v)^2\le \sum_{i=1}^N\sum_{v\in V} w_i(v)^2 \le NC.
\end{equation}
Now let $v_1\sim \dots \sim v_l$ be a chain connecting 
$V_E$ and $V_F$.
For $i\in \{1,\dots, N\}$ let 
 $m_i$ be the largest  index with $v_{m_i}\in V_{B_i}$.
Since $v_1\in V_F\sub V_{B_i}$ the number $m_i$ is well defined.
Moreover, $m_i \le m_{i+1}$. Let 
$m'_i$ be the smallest index $\ge  m_i$ with $v_{m'_i} \in V_{Z\setminus 
\la B_i}$. Note that $m'_i$ is well defined since  $v_l\in V_F
\sub V_{Z\setminus \la B_i}$.
Then $v_{m_i}\sim \dots \sim v_{m'_i}$  is a chain connecting 
$V_{B_i}$ and $V_{Z\setminus \la B_i}$ and we obtain from the 
admissibility of $w_i$ 
$$ \sum_{j=m_i}^{m'_i} w_i(v_j) \ge 1. $$ 

Let $i\in \{1,\dots, N-1\}$ and  let $j= m'_{i}$.
Assume  $m_i< m'_{i}$. Then 
 $v_{j-1} \notin V_{Z\setminus \la B_i}$ by definition of 
$m'_{i}$. This means 
$U_{v_{j-1}}\sub \la B_i $.
Then $ \emptyset \ne U_{v_{j-1}}\cap U_{v_{j}} \sub \la B_{i}\cap U_{v_{j}}$,
and so  $v_{j}\in V_{\la B_{i}}$. This is also true 
if $m'_{i}=m_i$.
 By  choice of $k$  we have 
$v_{j}\notin V_{Z\setminus B_{i+1}}$ which implies $ j< l$ and 
$U_{v_j}\sub  B_{i+1}$. Therefore, we have
that $\emptyset\ne U_{v_j}\cap U_{v_{j+1}}\sub B_{i+1} \cap U_{v_{j+1}}$.
Thus $v_{j+1}\in V_{B_{i+1}}$ and we conclude $m_{i+1} \ge j+1> m'_i$.  
In other words, the chains
$v_{m_i}\sim \dots \sim v_{m'_i}$ for $i\in \{1,\dots, N\}$ are pairwise
disjoint and we get 
$$ \sum_{\nu=1}^l w(v_\nu)\ge \sum_{i=1}^N \sum_{\nu =m_i}^{m'_{i}}
w_i(v_\nu)\ge N.$$    
We conclude that $w/N$ is admissible for the pair  $(V_E, V_F)$, 
and so by (\ref{mbd}) we have 
 $$\text{mod}_2(V_E, V_F)\le C/N.$$
Returning to the usual notation,  this means 
that  $\text{mod}_2^{G_k}(V^k_E, V^k_F)$ is small  for $k\ge k_3$,  
 if $\Delta(E,F)$
is large, quantitatively.

Proposition \ref{qs2} now shows that there exists an
$\tilde\eta$-quasisymmetric  map $f\: \sph \ra Z$, where $\tilde\eta$ depends
only on the data. This map has to be a homeomorphism. Its inverse
map will be  an $\eta$-quasisymmetric  homeomorphism $g\: Z\ra \sph$,
where $\eta$ depends only on the data. 

Conversely, suppose that $Z$ is quasisymmetric to $\sph$.
Assume that $Z$ is $\la_0$-$LLC$, where $\la_0> 1$.
By Theorem
\ref{necsuff} condition (\ref{uppbdmod}) will be satisfied for $L>0$
and  a suitable function $\Psi$. We can find $t_0> 0 $ and $C>0$ such that
$\Psi(t) <C$ for $t\ge t_0$. Let $\la:=2t_0+\la_0^2>1$. 
Suppose $B=B(a,r)$
is a ball in $Z$. From  $\la_0$-$LLC_1$ it follows that there exists
a  continuum $E$ with  $B\sub E\sub \bar B(a, \la_0 r)$.  Moreover, assume that 
$\la B\ne \emptyset$. Then $\la_0$-$LLC_2$ implies that there exists
a continuum $F$ with  $ Z\setminus \la B  \sub F \sub   Z\setminus
B(a, \la r/ \la_0)$. 
We have $\Delta(E,F)\ge (\la -\la_0^2)/2=t_0$. Since 
$\mesh({\cal A}_k)\to 0$, we have that $E$ and $F$ are not contained
in any $L$-star of ${\cal A}_k$  for sufficiently 
large $k$.
It follows that for these $k$ we have
$$ \GMod_2^{G_k}(V^k_{B}, V^k_{Z\setminus \la B}) \le  \GMod_2^{G_k}(V^k_E,
V^k_F)  < C.$$  If $Z\setminus \la B=\emptyset$, then 
$\GMod_2^{G_k}(V^k_{B}, V^k_{Z\setminus \la B})=0$ by definition of
the modulus. In any case we see that condition (\ref{uppbdmod2}) is satisfied. 
\qed

\subsection{Concluding remarks}
\label{concludingremarks}

\no(1) Theorems similar to Theorem \ref{mainthm} are   true for more 
general  surfaces. In  the case when $Z$ is homeomorphic to 
$\R^2$ the following statement holds:

{\em Let $Z$ be an Ahlfors $2$-regular metric space homeomorphic to 
$\R^2$. Then $Z$ is quasisymmetric to $\R^2$ (equipped with the 
standard Euclidean metric) if and only if $Z$ is   proper 
and  linearly locally connected.} 

Recall that a metric space is called {\em proper}
 if its closed balls are compact. 

\medskip\no
(2) Theorem \ref{mainthm} can be used to give a canonical model for 
$2$-regular $2$-spheres that are linearly locally contractible. 
To make this precise we remind the reader of the concept of a deformation 
of a metric  space $(Z,d)$  by a metric doubling measure. 
Suppose $\mu$ is a  Borel measure on $Z$. The measure is called 
{\em doubling} if there exists
a constant $C\ge 1$ such that 
$$ \mu(B(a,2r)\le C \mu(B(a,r)),$$
whenever $a\in Z$  and $r>0$.
If $x,y\in Z$ let $B_{xy}:=B(x, d(x,y))\cup B(y, d(x,y)).$  Suppose 
$Q\ge 1$ is fixed. Then we  introduce a  function  
$\delta_\mu(x,y):= \mu(B_{x,y})^{1/Q}$. The measure $\mu$ is called 
a {\em metric doubling measure} (with exponent $Q$)
if  $\delta_\mu$ is a metric up to 
a bounded multiplicative constant, i.e., there exists a metric $\delta$ 
on $Z$ and a constant $C\ge 1$ such that 
$$ (1/C) \delta(x,y)\le \delta_\mu(x,y) \le C \delta(x,y) \for x,y\in Z.$$ 
Suppose  $\mu$ is a  metric doubling measure. As long as an ambiguity 
caused by a multiplicative constant is harmless, the distance function
$\delta_\mu$ is as good as  a metric and we can talk about the metric 
space $(Z, \delta_\mu)$ and quasisymmetric maps of this space etc.
It is easy to see that  the 
``metric space" $(Z, \delta_\mu)$ is Ahlfors $Q$-regular and
quasisymmetric to $(Z,d)$ by the identity map. 

If $Z={\mathbb S}^n$  and $Q=n$, then every  metric doubling measure $\mu$
 is absolutely 
continuous with respect to spherical  measure  $\sigma_n$, i.e., there exists 
a measurable weight $w\: {\mathbb S}^n\to [0,\infty]$ such that 
$d\mu=w\, d\sigma_n$. The weight is an $A_\infty$-weight. Weights 
that arise from  metric doubling measures are called strong
$A_\infty$-weights.     

Theorem \ref{mainthm} now implies the following statement:

{\em A metric $2$-sphere  $(Z,d)$ is Ahlfors $2$-regular and linearly locally
contractible if and only if $(Z,d)$ is bilipschitz
to a space $(\sph, \delta_\mu)$, where $\mu$ is a metric doubling
measure on $\sph$ with exponent $Q=2$.}

Indeed, if $(Z,d)$ is Ahlfors $2$-regular and linearly locally
contractible, then there exists a quasisymmetric  homeomorphism
$f\:  \sph \ra Z$ by Theorem \ref{mainthm}. 
Define the   measure $\mu$ on $\sph$ as  the pull-back of ${\cal H}^2$ 
by $f$. So $\mu(E)={\cal H}^2(f(E))$ for a Borel set $E\sub \sph$.
Using the fact that $f$ is  quasisymmetric and that $Z$ is $2$-regular,
it easy to see that $\mu$ is doubling. Moreover, we have 
$\delta_\mu(x,y) \simeq d(f(x), f(y))$ for $x,y\in \sph$. This shows that 
$\mu$ is a metric doubling measure, and that $f\: (\sph, \delta_\mu) 
\ra (Z,d)$ is bilipschitz. 

Conversely, if $\mu$ is a metric doubling
measure on $\sph$ with exponent $Q=2$, then $(\sph,\delta_\mu)$ is $2$-regular.
Hence $(Z,d)$ is also $2$-regular, because this property is  preserved 
 under bilipschitz maps. Since $(Z,d)$ is bilipschitz to $(\sph,\delta_\mu)$
and the latter space is quasisymmetric to $\sph$ by the identity map,
the spaces $(Z,d)$ and $\sph$ are quasisymmetric. Linear local contractibility
is invariant under quasisymmetries, and since  $\sph$ has this property, 
so does $(Z,d)$.  

\medskip\no

(3) A necessary condition for a metric $2$-sphere  $Z$ to be bilipschitz
 to $\sph$ is that $Z$ is $2$-regular and linear locally 
contractible. By the   result in (3)  a space satisfying these
necessary conditions is bilipschitz to a space
$(\sph, \delta_\mu)$, where $\mu$ is a metric doubling measure on 
$\sph$  with exponent
$2$. So the problem of characterizing $\sph$ up to bilipschitz equivalence
is reduced to the question which of the spaces $(\sph, \delta_\mu)$ are 
bilipschitz  to $\sph$.   

This question is related to the Jacobian problem for quasiconformal 
mappings on $\sph$ as follows. If  $f\: \sph \ra \sph$ is a quasiconformal map, 
we  denote by $J_f$ its Jacobian (determinant).  The Jacobian
problem for quasiconformal  maps asks for a characterization 
of the weights
  $w\: \sph\to [0,\infty]$ for which there exists a
quasiconformal map $f\: \sph \ra \sph$ such that 
$$ (1/C) J_f(x) \le w(x)\le C J_f(x) \for \sigma_2\text{-a.e. } x\in \sph,$$
where $C$ is a constant independent of $x$. 
A necessary and sufficient  condition for a weight $w$ to be comparable to 
a Jacobian of  a quasiconformal map is that $w$ is a strong $A_\infty$-weight,
i.e., the measure $\mu$ defined by $d\mu=w\, d\sigma_2$ is a metric 
doubling measure, and that $(\sph, \delta_\mu)$ is bilipschitzly 
equivalent to $\sph$.

{}From this  we see  that the  Jacobian problem for quasiconformal mappings on  
$\sph$ is equivalent with  the problem of characterizing 
$\sph$ up to bilipschitz equivalence.

\medskip\no 
(4) The usefulness of Theorem \ref{suff} 
depends on whether one can verify its hypotheses in 
concrete situations. There are some interesting fractal spaces of Hausdorff 
dimension greater than $2$ where this can be done. 
For example, consider the space $Z\sub \R^3$ 
obtained as follows. The space $Z$ will be the limit of a sequence 
of two-dimensional  cell complexes $Z_n$. Each $Z_n$ 
consists  of a union 
of congruent oriented squares. The orientation of each square is visualized 
by specifying which of the two  directions perpendicular  to the square
is considered as normal.    
The sets $Z_n$ are inductively constructed as follows. 
The cell complex 
$Z_0$ is the boundary of the unit cube $I^3\sub\R^3$, 
 where the $2$-cells 
are the  six squares forming the faces of $Z_0$. We orient the 
squares of $Z_0$ by assigning to them the normal pointing outward 
$I^3$.
Now $Z_{n+1}$ is obtained from $Z_n$  by modifying 
each of the oriented  squares $S$ forming $Z_n$ as follows. 
Subdivide $S$ into 25 congruent subsquares with the induced 
orientation. (Actually any fixed number
$(2k+1)^2$ with $k\ge 2$ could be taken here. In the case $k=1$ there 
are some problems with overlaps in the inductive construction.)  
On the  ``central" subsquare $S'$
of $S$
 place an appropriately sized  cube $C$ in the normal direction
so that one of the faces of $C$ agrees with $S'$. 
The    face squares of $C$ are   oriented so that 
their normals point outward $C$.  The  desired 
modification of $S$ is now obtained 
by replacing  the  ``central" subsquare $S'$ of $S$ by the oriented faces
of $C$ different from $S'$ and keeping all other oriented subsquares. In this 
way each square of $Z_n$ leads to $24+5=29$ squares of $Z_{n+1}$.
The limit set $Z$ is equipped with the ambient metric of $\R^3$.
It can be shown that $Z$ is homeomorphic to $\sph$ and $Q$-regular 
for some $Q>2$.  Using the symmetry properties of $Z$ and   Theorem 
\ref{suff},  one can show:
{\em $Z$ is quasisymmetric to $\sph$}.
An independent proof of this fact based on the dynamics of rational functions
is due to D.\ Meyer \cite{Meyer}.

 We hope to explore applications of 
Theorem \ref{suff} more systematically in the future.

\bibliography{refs2}
\bibliographystyle{siam}
\addcontentsline{toc}{subsection}{References}

\end{document}